\documentclass[reqno,11pt]{amsart}
 
\usepackage{amsmath,amsthm}
\usepackage{amssymb}
\usepackage{euscript}

\numberwithin{equation}{subsection}

\newcommand{\sqsp}{\renewcommand{\baselinestretch}{1.08}\tiny\normalsize}

\newtheorem{thm}[subsection]{Theorem}
\newtheorem{lemma}[subsection]{Lemma}
\newtheorem{prop}[subsection]{Proposition}
\newtheorem{cor}[subsection]{Corollary}
\newtheorem{definition}[subsection]{Definition}

\newcommand{\cat}[1]{{\EuScript #1}}
\newcommand{\cF}{\cat{F}}
\newcommand{\cC}{\cat{C}}
\newcommand{\cP}{\cat{P}}

\newcommand{\bZ}{\mathbf{Z}}
\newcommand{\bN}{\mathbf{N}}
\newcommand{\Hl}{H_\lambda}
\newcommand{\Hbar}{\bar{H}_\lambda}
\newcommand{\Der}{\mathbf{Der}_\lambda}
\newcommand{\Derbar}{\overline{\mathbf{Der}}_\lambda}
\newcommand{\Inn}{\mathbf{Inn}_\lambda}
\newcommand{\Innbar}{\overline{\mathbf{Inn}}_\lambda}

\DeclareMathOperator{\Id}{Id} 
\DeclareMathOperator{\End}{\mathbf{End}}
\DeclareMathOperator{\Endbar}{\overline{\mathbf{End}}}

\begin{document}
\title{Cohomology of truncated polynomial $\lambda$-rings}
\author{Donald Yau}

\begin{abstract}
The $\lambda$-ring cohomology, $\Hl^0$ and $\Hl^1$, of certain truncated polynomial filtered $\lambda$-rings are computed.
\end{abstract}

\email{dyau@math.uiuc.edu}
\address{Department of Mathematics, University of Illinois at Urbana-Champaign, 1409 W. Green Street, Urbana, IL 61801 USA}
\date{\today}

\maketitle

\sqsp


\section{Introduction}
\label{sec:intro}

Lambda-rings were introduced by Grothendieck \cite{gro} to describe algebraic objects equipped with operations $\lambda^i$ that behave like exterior powers.  Since their introduction by Grothendieck, $\lambda$-rings have been shown to play important roles in several areas of mathematics, including Topology, Representation Theory, and Algebra.  Indeed, examples of $\lambda$-rings include the complex $K$-theory of a topological space, the representation ring of a group, and the universal Witt ring of a commutative ring.  See the references \cite{at,gro,hazewinkel,knutson,patras}.  Note that what we call a $\lambda$-ring here is sometimes referred to as a \emph{special} $\lambda$-ring in the literature.

Any reasonable class of algebras (associative, Lie, Poisson, Jordan, von Neumann, etc.) should come with (i) a cohomology theory $H^*$ and (ii) an algebraic deformation theory which is describable in terms of $H^*$.  The original theory of algebraic deformations for associative algebras was worked out by Gerstenhaber \cite{ger}, and the relevant cohomology theory is Hochschild cohomology.  Deformations and the associated cohomology for many other types of algebras have since been worked out.  Since $\lambda$-rings form such a useful class of algebras, there should likewise be deformations and cohomology for $\lambda$-rings.  Indeed, the author defined in \cite{yau3} $\lambda$-ring cohomology $\Hl^*(R)$ for a $\lambda$-ring $R$ and used it to study $\lambda$-ring deformations along the path initiated by Gerstenhaber \cite{ger}.  There is even a ``composition" product on the cochain complex which gives rise to $\lambda$-ring cohomology, and this product descends to $\Hl^*(R)$, giving it the structure of a graded, associative, unital algebra.

The current paper can be considered as a sequel to \cite{yau3}.  The main purpose here is to show how $\lambda$-ring cohomology, at least in dimensions $0$ and $1$, can be computed for some filtered $\lambda$-ring structures on truncated polynomial algebras.  Particular attention is paid to those filtered $\lambda$-ring structures that are candidates for topological realizations (i.e.\ those that can possibly be the $K$-theory of a torsionfree space).

Low-dimensional $\lambda$-ring cohomology is intrinsically interesting.  In fact, under the composition product, $\Hl^0(R)$ is a non-trivial subalgebra of the algebra of $\bZ$-linear self maps of $R$ that commute with all the Adams operations.  Each Adams operation of $R$ does lie in $\Hl^0(R)$.  The group $\Hl^1(R)$ is important because it is the natural habitat for the ``infinitesimal" (sort of like the ``initial velocity") of a $\lambda$-ring deformation of $R$.  Indeed, $\Hl^1(R)$ classifies the order $1$ $\lambda$-ring deformations of $R$ up to equivalence.

A description of the rest of the paper follows.

Section \ref{sec:review} gives a brief account of the basics of (filtered) $\lambda$-rings and Adams operations.  It also reviews $\lambda$-ring cohomology and the composition product.

The main result in Section \ref{sec:H1} is Theorem \ref{thm:H1group}, which shows that $\Hl^1$ is isomorphic to a smaller group $\Hbar^1$ that involves only the Adams operations $\psi^p$ for $p$ primes (as opposed to all the $\psi^k$).  This will allow a more efficient computation of the group $\Hl^1$ in later sections.  This section ends with a description of the graded algebra $\Hl^{\leq 1}(\bZ) = \Hl^*(\bZ)/(\oplus_{i \geq 2}\, \Hl^i(\bZ))$ for the $\lambda$-ring $\bZ$ of integers.

Section \ref{sec:n=2} gives a description of the graded algebra $\Hl^{\leq 1}(R)$ for each of the uncountably many isomorphism classes of filtered $\lambda$-ring structures $R$ on the dual number ring $\bZ \lbrack x \rbrack/(x^2)$ (Theorem \ref{thm:dual}).

Section \ref{sec:H0 n=3} begins by recalling the uncountably many isomorphism classes of filtered $\lambda$-ring structures on the truncated polynomial ring $\bZ \lbrack x \rbrack/(x^3)$.  The main result of this section is Theorem \ref{thm:H0 n=3}, which computes the algebra $\Hl^0(R)$ for each one of these isomorphism classes of filtered $\lambda$-rings.  It turns out that each such $\Hl^0(R)$ is a $3$-dimensional commutative subalgebra of the algebra $M(3, \bZ)$ of $3$-by-$3$ integer matrices.

Among the uncountably many isomorphism classes of filtered $\lambda$-ring structures on the truncated polynomial ring $\bZ \lbrack x \rbrack/(x^3)$, at most $64$ of them can be realized as the $K$-theory of torsionfree spaces (see \cite{yau4}).  The main result of section \ref{sec:H1 n=3} is Theorem \ref{thm:H1 n=3}, which describes the groups $\Hl^1$ and determines (non-)commutativity of the graded algebras $\Hl^{\leq 1}$ for these $64$ isomorphism classes of filtered $\lambda$-rings.  Interestingly enough, exactly $35$ of those $64$ graded algebras $\Hl^{\leq 1}(R)$ are commutative.

It is also shown in \cite{yau4} that the truncated polynomial ring $\bZ \lbrack x \rbrack/(x^4)$ admits uncountably many isomorphism classes of filtered $\lambda$-ring structures and that, among these classes, at most $61$ of them can be realized by the $K$-theory of torsionfree spaces.  The main results of section \ref{sec:n=4}, Theorem \ref{thm1:H0 n=4} and Theorem \ref{thm2:H0 n=4}, compute the algebras $\Hl^0$ for these $61$ isomorphism classes of filtered $\lambda$-rings.  Each such $\Hl^0$ turns out to be a $4$-dimensional commutative subalgebra of the algebra $M(4, \bZ)$ of $4$-by-$4$ matrices with integer entries.


\section{$\lambda$-rings and cohomology}
\label{sec:review}

The purpose of this section is to review some basic definitions about (filtered) $\lambda$-rings, Adams operations, and $\lambda$-ring cohomology.  For more discussions about basic properties of $\lambda$-rings, consult the references \cite{at,knutson,patras}.  The reference for $\lambda$-ring cohomology is the author's paper \cite{yau3}.

\subsection{$\lambda$-rings}
\label{subsec:lambda}

By a $\lambda$-ring we mean a unital, commutative ring $R$ equipped with functions 
   \[
   \lambda^i \colon R ~\to~ R \quad (i \geq 0),
   \]
called $\lambda$-operations, which satisfy the following conditions.  For any integers $i, j \geq 0$ and elements $r$ and $s$ in $R$:
   \begin{itemize}
   \item $\lambda^0(r) = 1$.
   \item $\lambda^1(r) = r$.
   \item $\lambda^i(1) = 0$ for $i > 1$.
   \item $\lambda^i(r + s) = \sum_{k = 0}^i\, \lambda^k(r)\lambda^{i-k}(s)$.
   \item $\lambda^i(rs) = P_i(\lambda^1(r), \ldots , \lambda^i(r); \lambda^1(s), \ldots, \lambda^i(s))$.
   \item $\lambda^i(\lambda^j(r)) = P_{i,j}(\lambda^1(r), \ldots , \lambda^{ij}(r))$.
   \end{itemize}
The $P_i$ and $P_{i,j}$ are some universal polynomials with integer coefficients, defined as follows.  Consider the variables $\xi_1, \ldots , \xi_i$ and $\eta_1, \ldots , \eta_i$.  Denote by $s_1, \ldots , s_i$ and $\sigma_1, \ldots , \sigma_i$, respectively, the elementary symmetric functions of the $\xi$'s and the $\eta$'s.  The polynomial $P_i$ in $2i$ variables is defined by insisting that the expression 
$P_i(s_1, \ldots , s_i; \sigma_1, \ldots , \sigma_i)$ 
be the coefficient of $t^i$ in the finite product
   \[
   \prod_{m,n=1}^i\, (1 + \xi_m \eta_n t).
   \]
Similarly, if $s_1, \ldots , s_{ij}$ are the elementary symmetric functions of $\xi_1, \ldots , \xi_{ij}$, then the polynomial $P_{i,j}$ in $ij$ variables is defined by insisting that the expression 
$P_{i,j}(s_1, \ldots , s_{ij})$ 
be the coefficient of $t^i$ in the finite product
   \[
   \prod_{l_1 < \cdots < l_j} \, (1 + \xi_{l_1} \cdots \xi_{l_j} t).
   \]

Note that what we call a $\lambda$-ring here is referred to as a \emph{special} $\lambda$-ring in \cite{at,gro}.

A map of $\lambda$-rings is a ring map that commutes with the $\lambda$-operations.

\subsubsection{Adams operations}
\label{subsec:adams}

In a $\lambda$-ring $R$, there are the so-called Adams operations
   \[
   \psi^n \colon R ~\to~ R \quad (n \geq 1),
   \]
satisfying the following properties:
\begin{itemize}
\label{adams}
\item All the $\psi^n$ are ring maps.
\item $\psi^1 = \Id$.
\item $\psi^m \psi^n = \psi^{mn} = \psi^n \psi^m$.
\item $\psi^p(r) \equiv r^p$ (mod $pR$) for each prime $p$ and element $r$ in $R$.
\end{itemize}
The Adams operations can be defined inductively by the Newton formula:
   \[
   \psi^n(r) - \lambda^1(r)\psi^{n-1}(r) + \cdots + (-1)^{n-1}\lambda^{n-1}(r) \psi^1(r) + (-1)^n n\lambda^n(r) = 0.
   \]

Suppose given a unital, commutative, $\bZ$-torsionfree ring $R$ with self ring maps $\psi^n \colon R \to R$ satisfying the above four properties of Adams operations.  Then a result of Wilkerson \cite{wil} says that there exists a unique $\lambda$-ring structure on $R$ whose Adams operations are exactly the given $\psi^n$.

\subsubsection{Filtered $\lambda$-rings}

By a \emph{filtered ring} we mean a (unital, commutative) ring $R$ equipped with a decreasing sequence of ideals $I^n$.  A map of filtered rings is a ring map that preserves the filtration ideals.

A \emph{filtered} $\lambda$-\emph{ring} is a filtered ring $R$ that is also a $\lambda$-ring in which the filtration ideals are all closed under the $\lambda$-operations, i.e.\ $\lambda^i(I^n) \subseteq I^n$ for all $n$ and all $i > 0$.  A
 map of filtered $\lambda$-rings is a $\lambda$-ring map that is also a filtered ring map.

For example, the $K$-theory of a topological space $X$ with the homotopy type of a CW complex is a filtered $\lambda$-ring.  Here the filtration ideals are given by the kernels of the restriction maps to the skeletons, i.e.\ 
   \[
   I^n ~=~ \ker(K(X) ~\to~ K(X_{n-1})),
   \]
where $X_{n-1}$ is the $(n-1)$-skeleton of $X$ and the map is induced by the inclusion $X_{n-1} \subset X$.  The Cellular Approximation Theorem assures that any two CW structures on $X$ give rise to isomorphic filtered $\lambda$-ring structures on $K(X)$.

The obvious analogue of Wilkerson's Theorem discussed in section \ref{subsec:adams} holds for the truncated polynomial algebras $\bZ \lbrack x \rbrack/(x^n)$ with $x$ in some fixed positive filtration.

\subsection{$\lambda$-ring cohomology}
\label{subsec:coh}

Let $R$ be a $\lambda$-ring with $\lambda$-operations $\lambda^i$ $(i \geq 0)$ and Adams operations $\psi^n$ $(n \geq 1)$.

\subsubsection{The complex $\cF^*$}
\label{subsubsec:F}

The $\lambda$-ring cohomology groups $\Hl^*(R)$ of $R$ are defined to be the cohomology groups of a certain cochain complex $\cF^*(R)$, which is defined as follows.

Denote by $\End(R)$ the algebra of $\mathbf{Z}$-linear endomorphisms of $R$. Let $\Endbar(R)$ be the subalgebra of $\End(R)$ consisting of those linear endomorphisms $f$ of $R$ that satisfy the condition,
   \[
   \label{eq:endbar def}
   f(r)^p ~\equiv~ f(r^p) \pmod{pR},
   \]
for every prime $p$ and each element $r \in R$.  Note that every self ring map of $R$ lies in $\Endbar(R)$.  In particular, $\Endbar(R)$ contains all the Adams operations of $R$.

We are now ready to define the cochain complex $\cF^* = \cF^*(R)$.  Let $T$ be the set of positive integers.  Define $\cF^0$ to be the underlying additive group of $\Endbar(R)$ and $\cF^1$ to be the set of functions
   \[
   \label{eq:F1}
   f \colon T ~\to~ \End(R)
   \]
satisfying the condition
   \[
   f(p)(R) ~\subseteq~ pR 
   \]
for every prime $p$.  For $n \geq 2$, $\cF^n$ is defined to be the set of functions
   \[
   f \colon T^n ~\to~ \End(R).
   \]
If $f$ and $g$ are elements of $\cF^n$, then their sum is defined by 
   \[
   (f + g)(m_1, \ldots, m_n)(r) ~=~ f(m_1, \ldots, m_n)(r) + g(m_1, \ldots, m_n)(r)
   \]
for $(m_1, \ldots, m_n) \in T^n$ and $r \in R$.  This gives $\cF^n$ $(n \geq 0)$ a natural additive group structure.

The differential 
   \[
   d^n \colon \cF^n ~\to~ \cF^{n+1}
   \]
is defined by the formula
   \begin{multline*}
   \label{eq:dn}
   (d^n f)(m_0, \ldots, m_n)  \\
   ~=~ \psi^{m_0} f(m_1, \ldots, m_n) 
   ~+~ \sum_{i=1}^n \, (-1)^i f(m_0, \ldots, m_{i-1}m_i, \ldots, m_n) \\
   ~+~ (-1)^{n+1} f(m_0, \ldots, m_{n-1}) \psi^{m_n}.
   \end{multline*}
The $d^n$ are linear maps and satisfy $d^{n+1} d^n = 0$ for each $n \geq 0$.  This makes $\cF^* = (\cF^*(R), d^*)$ into a cochain complex.  We define the $n$th $\lambda$-\emph{ring cohomology group of} $R$, denoted by $H^n_\lambda(R)$, to be the $n$th cohomology group of the cochain complex $\cF^* = (\cF^*(R), d^*)$.

\subsubsection{Composition product}

Given a $\lambda$-ring $R$, there is an associative, bilinear pairing
   \[
   - \circ -  ~\colon~ \cF^n \otimes \cF^k ~\to~ \cF^{n+k} \quad (n, k \geq 0)
   \]
on the complex $\cF^* = (\cF^*(R), d^*)$ defined by
   \[
   (f \circ g)(m_1, \ldots, m_{n+k}) 
   ~=~ f(m_1, \ldots, m_n) \circ g(m_{n+1}, \ldots, m_{n+k}),
   \]
for $f \in \cF^n$ and $g \in \cF^k$.  The element $\Id_R \in \cF^0$ acts as a two-sided identity for this pairing.  Moreover, this pairing satisfies the Leibnitz identity,
   \[
   d(f \circ g) ~=~ (df) \circ g ~+~ (-1)^{\vert f \vert}f \circ (dg),
   \]
where $\vert f \vert$ is the dimension of $f$.  We call this pairing the \emph{composition product}.  The complex $\cF^*$ with the composition product is a differential graded algebra.

The Leibnitz identity implies that the composition product descends to $H^*_\lambda(R) = \oplus_i \, \Hl^i(R)$ with 
   \[
   \lbrack f \rbrack \circ \lbrack g \rbrack 
   ~=~ \lbrack f \circ g \rbrack, 
   \]
where $\lbrack f \rbrack$ denotes the cohomology class of a cocycle.  With this product, the graded group $\Hl^*(R)$ becomes a graded, associative, unital algebra.


\section{A reinterpretation of $\Hl^1$}
\label{sec:H1}

Since $\lambda$-ring cohomology are defined using the Adams operations, which are determined by $\psi^p$ for $p$ primes, it should be the case that $\lambda$-ring cohomology can be reinterpreted in terms of only the $\psi^p$.  This main point of this section is to do exactly that for $\Hl^1$.

Let $\cP$ denote the set of all primes.  For a $\lambda$-ring $R$, define the following:

\begin{itemize}
\item $\Der(R) =$ The set of $f \in \cF^1(R)$ satisfying the condition
   \begin{equation}
   \label{eq:f(mn)} 
   f(mn) ~=~ \psi^m f(n) ~+~ f(m) \psi^n
   \end{equation}
for all positive integers $m$ and $n$.  This is called the group of $\lambda$-\emph{derivations} in $R$.  Note that $f(1)$ must be $0$.

\item $\Derbar(R) =$ The set of sequences $\lbrace f(p) \rbrace_{p \in \cP}$ indexed by the primes with each $f(p) \in \End(R)$, such that $f(p)(R) \subseteq pR$ and that
   \begin{equation}
   \label{eq:f(pq)}
   \psi^p  f(q) ~+~ f(p) \psi^q ~=~ \psi^q f(p) ~+~ f(q) \psi^p
   \end{equation}
for all primes $p$ and $q$.  This forms a group under coordinatewise addition.

\item $\Inn(R) =$ The set of elements in $\cF^1(R)$ of the form $\lbrack \psi^*, g \rbrack$, where $g \in \cF^0(R)$ and  $\lbrack \psi^*, g \rbrack(n) = \psi^n g - g \psi^n$.  This is called the group of $\lambda$-\emph{inner derivations} in $R$.

\item $\Innbar(R) =$ The set of sequences $\lbrace \lbrack \psi^p, g \rbrack \rbrace_{p \in \cP}$ indexed by the primes in which $g \in \cF^0(R)$.  This is clearly a subgroup of $\Derbar(R)$.

\item $\Hbar^1(R) = $ The quotient group $\Derbar(R)/\Innbar(R)$.
\end{itemize}

It follows directly from the definition that $\Der(R)$ is the kernel of $d^1$ and that $\Inn(R)$ is the image of $d^0$.  Therefore, we have

\begin{prop}[= Proposition 9 in \cite{yau3}]
As an additive group, $\Hl^1(R)$ is the quotient $\Der(R)/\Inn(R)$.
\end{prop}

Notice that \eqref{eq:f(mn)} implies \eqref{eq:f(pq)}, and so there is a linear map
   \[
   \pi \colon \Der(R) ~\to~ \Derbar(R), 
   \]
where
   \[
   \pi(f) ~=~ \lbrace f(p) \rbrace_{p \in \cP}.
   \]
The image of $\Inn(R)$ under $\pi$ is exactly $\Innbar(R)$.  It follows that $\pi$ induces a linear map
   \begin{equation}
   \label{eq:pi}
   \pi \colon \Hl^1(R) ~\to~ \Hbar^1(R)
   \end{equation}
which sends $\lbrack f \rbrack$ to $\lbrack \lbrace f(p) \rbrace_{p \in \cP} \rbrack$.

\begin{thm}
\label{thm:H1group}
The map $\pi$ in \eqref{eq:pi} is a group isomorphism.
\end{thm}

From now on, using the isomorphism $\pi$, we will consistently identify $\Hl^1(R)$ with $\Hbar^1(R)$.

\begin{proof}
Theorem \ref{thm:H1group} will follow from Lemma \ref{lem1:H1group} and Lemma \ref{lem2:H1group} below and the fact that $\pi(\Inn(R)) = \Innbar(R)$.
\end{proof}

\begin{lemma}
\label{lem1:H1group}
Let $f$ be an element of $\Der(R)$ and let $n = p_1 \cdots p_r$ be a positive integer, in which the $p_i$ are primes, not necessarily distinct.  Then we have that
   \begin{equation}
   \label{eq:f(n)}
   f(n) ~=~ \sum_{i=1}^r\, \psi^{p_1 \dots p_{i-1}}  f(p_i)  \psi^{p_{i+1}\cdots p_r}.
   \end{equation}
In particular, the map $\pi \colon \Der(R) \to \Derbar(R)$ is injective.
\end{lemma}

\begin{proof}
We proceed by induction on the number $r$ of prime factors.  There is nothing to prove when $r = 1$.  Suppose that \eqref{eq:f(n)} has been proved for positive integers $< r$.  Write $m$ for $n/p_r = p_1 \cdots p_{r-1}$.  Then we have
   \[
   \begin{split}
   f(n) &~=~ f(mp_r) \\
        &~=~ \psi^m f(p_r) + f(m) \psi^{p_r} \\
        &~=~ \left(\sum_{i=1}^{r-1}\, \psi^{p_1 \cdots p_{i-1}} f(p_i) \psi^{p_{i+1}\cdots p_{r-1}}\right) \psi^{p_r} ~+~ \psi^m f(p_r) \\
        &~=~ \sum_{i=1}^r \, \psi^{p_1 \cdots p_{i-1}} f(p_i) \psi^{p_{i+1}\cdots p_r}.
   \end{split}
   \]
This finishes the induction and proves the Lemma.
\end{proof}

\begin{lemma}
\label{lem2:H1group}
Given an element $\lbrace g(p) \rbrace_p \in \Derbar(R)$, there exists an element $f \in \Der(R)$ such that $f(p) = g(p)$ for each prime $p$.  In other words, the map $\pi \colon \Der(R) \to \Derbar(R)$ is surjective.
\end{lemma}

From now on, using the isomorphism $\pi$ and Lemma \ref{lem1:H1group} and Lemma \ref{lem2:H1group}, we will consistently identify $\Der(R)$ (respectively $\Inn(R)$) with $\Derbar(R)$ (respectively $\Innbar(R)$).

\begin{proof}
Set $f(1) = 0$.  Given an integer $n = p_1 \cdots p_r > 1$, where each $p_i$ is a prime, define $f(n)$ by
   \begin{equation}
   \label{eq:g(p)}
   f(n) ~\buildrel \text{def} \over =~ \sum_{i=1}^r \, \psi^{p_1\cdots p_{i-1}} g(p_i) \psi^{p_{i+1} \cdots p_r}.
   \end{equation}
We first have to make sure that this is well-defined.  In other words, we have to show that the right-hand side of \eqref{eq:g(p)} is independent of the order of the $p_i$ appearing in the prime factorization of $n$.  This is clearly true when $r = 1$; the case $r = 2$ follows from \eqref{eq:f(pq)}.

Suppose that $r \geq 2$.  Pick any $j$ with $1 \leq j \leq n - 1$ and write 
   \[
   n ~=~ p_1 \cdots p_{j-1}p_{j+1}p_jp_{j+2} \cdots p_r,
   \]
i.e., transposes $p_j$ and $p_{j+1}$.  With this way of writing $n$, we have
   \begin{multline}
   \label{eq:ind}
   f(n) ~=~ \left(\sum_{i=1}^{j-1}\, \psi^{p_1 \cdots p_{i-1}} f(p_i) \psi^{p_{i+1} \cdots p_r}\right) ~+~ \psi^{p_1 \cdots p_{j-1}}f(p_{j+1})\psi^{p_jp_{j+2} \cdots p_r} \\ 
   ~+~ \psi^{p_1\cdots p_{j-1}p_{j+1}}f(p_j)\psi^{p_{j+2}\cdots p_r} ~+~ \left(\sum_{l=j+2}^{r} \, \psi^{p_1\cdots p_{l-1}}f(p_l)\psi^{p_{l+1} \cdots p_r}\right).
   \end{multline}
The sum of the two terms in the middle (those that are not surrounded by parentheses) is
   \[
   \begin{split}
   \psi^{p_1 \cdots p_{j-1}} f(p_{j+1})& \psi^{p_jp_{j+2} \cdots p_r} ~+~ \psi^{p_1\cdots p_{j-1}p_{j+1}} f(p_j) \psi^{p_{j+2}\cdots p_r} \\
   &~=~ \psi^{p_1 \cdots p_{j-1}}\left(\psi^{p_{j+1}}f(p_j) ~+~ f(p_{j+1}) \psi^{p_j}\right)\psi^{p_{j+2} \cdots p_r} \\
   &~=~ \psi^{p_1 \cdots p_{j-1}}\left(\psi^{p_{j+1}}g(p_j) ~+~ g(p_{j+1})\psi^{p_j}\right)\psi^{p_{j+2} \cdots p_r} \\
   &~=~ \psi^{p_1 \cdots p_{j-1}}\left(\psi^{p_j}g(p_{j+1}) ~+~ g(p_j)\psi^{p_{j+1}}\right)\psi^{p_{j+2} \cdots p_r} \\
   &~=~ \psi^{p_1 \cdots p_{j-1}}\left(\psi^{p_j}f(p_{j+1}) ~+~ f(p_j)\psi^{p_{j+1}}\right)\psi^{p_{j+2} \cdots p_r}. \\
   \end{split}
   \]
Therefore, the right-hand sides of \eqref{eq:g(p)} and \eqref{eq:ind} agree.  Since every permutation on the set $\lbrace 1, \ldots , r \rbrace$ can be written as a product of transpositions of the form $(j, j+1)$ for $1 \leq j \leq r - 1$, the argument above shows that \eqref{eq:g(p)} is indeed well-defined.

Since $f(p) = g(p)$ for all primes $p$, to show that $f = \lbrace f(n) \rbrace$ lies in $\Der(R)$, it remains to show that $f$ satisfies \eqref{eq:f(mn)}.  Given $m = p_1 \cdots p_r$ and $n = q_1 \cdots q_s$, where the $p_i$ and $q_j$ are primes, we have
   \[
   \begin{split}
   &f(mn) ~=~ f(p_1 \cdots p_r q_1 \cdots q_s) \\
   &~=~ \sum_{i=1}^r \, \psi^{p_1 \cdots p_{i-1}} f(p_i) \psi^{p_{i+1} \cdots p_r q_1 \cdots q_s} ~+~ \sum_{j=1}^s \, \psi^{p_1 \cdots p_rq_1 \cdots q_{j-1}} f(q_j) \psi^{q_{j+1} \cdots q_s} \\
   &~=~ \left(\sum_{i=1}^r\, \psi^{p_1 \cdots p_{i-1}} f(p_i) \psi^{p_{i+1} \cdots p_r}\right) \psi^n ~+~ \psi^m\left(\sum_{j=1}^s \, \psi^{q_1 \cdots q_{j-1}} f(q_j) \psi^{q_{j+1} \cdots q_s}\right) \\
   &~=~ \psi^mf(n) ~+~ f(m)\psi^n.
   \end{split}
   \]
This finishes the proof of the Lemma.
\end{proof}

Consider the quotient
   \[
   \Hl^{\leq 1}(R) ~=~ \Hl^0(R) \oplus \Hl^1(R) ~\cong~ \frac{\Hl^*(R)}{\oplus_{i=2}^\infty\, \Hl^i(R)}
   \]
of the graded algebra $\Hl^*(R)$ (under the composition product) by the homogeneous ideal of elements of degree at least $2$.  We still consider $\Hl^{\leq 1}(R)$ a graded algebra, with $\Hl^0(R)$ and $\Hl^1(R)$ in degrees $0$ and $1$, respectively, and $0$ in degrees $\not= 0$, $1$.  Similarly, define the graded group
   \[
   \Hbar^{\leq 1}(R) ~\buildrel \text{def} \over =~ \Hl^0(R) \oplus \Hbar^1(R),
   \]
in which $\Hl^0(R)$ and $\Hbar^1(R)$ are in degrees $0$ and $1$, respectively.

With these notations, an immediate consequence of Theorem \ref{thm:H1group} is

\begin{cor}
\label{cor1:H1}
The group isomorphism
   \begin{equation}
   \label{eq:id+pi}
   (\Id,\, \pi) \colon \Hl^{\leq 1}(R) ~\xrightarrow{\cong}~ \Hbar^{\leq 1}(R)
   \end{equation}
induces a (unital, associative) graded algebra structure on $\Hbar^{\leq 1}(R)$ that is isomorphic to $\Hl^{\leq 1}(R)$.  Under this graded algebra structure of $\Hbar^{\leq 1}(R)$, one has that, for $g \in \Hl^0(R)$ and $\lbrace f(p) \rbrace_p \in \Derbar(R)$,
   \[
   g \circ \lbrack \, \lbrace f(p) \rbrace_{p\in \cP} \, \rbrack 
   ~=~ \lbrack \, \lbrace g f(p) \rbrace_{p \in \cP} \, \rbrack
   \]
and
   \[
    \lbrack \, \lbrace f(p) \rbrace_{p \in \cP} \, \rbrack \circ g
   ~=~  \lbrack \, \lbrace f(p) g \rbrace_{p \in \cP} \, \rbrack.
   \]
\end{cor}

Here $\circ$ denotes the algebra product in $\Hbar^{\leq 1}(R)$, the composition product, and $gf(p)$ is just composition of $\bZ$-linear self maps of $R$.

Using Corollary \ref{cor1:H1}, we will consistently identify the graded algebras $\Hl^{\leq 1}(R)$ and $\Hbar^{\leq 1}(R)$.

Recall that the ring $\bZ$ of integers has a unique $\lambda$-ring structure with $\lambda^i(n) = \binom{n}{i}$ and $\psi^k = \Id$ for all $k$.  Using Corollary \ref{cor1:H1}, the result in \cite{yau3} describing the groups $\Hl^{0}(\bZ)$ and $\Hl^1(\bZ)$ can be restated as follows.

\begin{thm}[= Corollary 8 and Corollary 10 in \cite{yau3}]
\label{thm:Z}
There is a graded algebra isomorphism 
   \[
   \Hl^{\leq 1}(\bZ) ~\cong~ \bZ \,\oplus\, \prod_{p \in \cP}\, p\bZ,
   \]
in which $\bZ$ is the degree $0$ part and the other summand is the degree $1$ part.  Given $m, n \in \bZ$ and $(n_p) \in \prod_{p \in \cP}\, p\bZ$, we have
   \[
   m \circ n ~=~mn
   \]
and
   \[
   m \circ (n_p) ~=~ (mn_p) ~=~ (n_p) \circ m.
   \]
In particular, $\Hl^{\leq 1}(\bZ)$ is a commutative graded algebra.
\end{thm}


\section{The dual number ring}
\label{sec:n=2}

Recall from \cite[Corollary 4.1.2]{yau1} that the dual number ring $\bZ \lbrack x \rbrack/(x^2)$, with $x$ in some fixed positive filtration $d$, admits uncountably many isomorphism classes of filtered $\lambda$-ring structures.  In fact, there is a bijection between this set of isomorphism classes and the set of sequences $(b_p)$ indexed by the primes in which $b_p \in p\bZ$.  The (isomorphism class of) filtered $\lambda$-ring $R$ corresponding to a sequence $(b_p)$ has Adams operations
   \[
   \psi^p(x) ~=~ b_px
   \]
for each prime $p$.

We will continue to denote by $\cP$ the set of all primes.

\begin{thm}
\label{thm:dual}
Let $R$ be a representative of any one of the uncountably many isomorphism classes of filtered $\lambda$-ring structures on the dual number ring $\bZ \lbrack x \rbrack/(x^2)$.  Then there is a graded algebra isomorphism 
   \[
   \Hl^{\leq 1}(R) ~\cong~ \bZ \oplus \bZ \,\oplus\, \prod_{p \in \cP}\, (p\bZ \oplus p\bZ),
   \]
in which $\bZ \oplus \bZ$ is the degree $0$ part and the infinite product is the degree $1$ part.  Given $(a, b),\, (c, d) \in \bZ \oplus \bZ$ and $\lbrace (x_p, y_p) \rbrace \in \prod_{p \in \cP}\, p\bZ \oplus p\bZ$, we have
   \[
   (a,\, b) \circ (c, \, d) ~=~ (ac, \, bd)
   \]
and
   \[
   (a,\, b) \circ \lbrace (x_p,\, y_p) \rbrace 
   ~=~ \lbrace (ax_p,\, by_p) \rbrace
   ~=~ \lbrace (x_p,\, y_p) \rbrace \circ (a,\, b).
   \]
In particular, $\Hl^{\leq 1}(R)$ is a commutative graded algebra.
\end{thm}

\begin{proof}
\emph{Step 1}: $\Hl^0(R)$.

We begin by computing $\Endbar(R)$, which, by definition, consists of the $\bZ$-linear maps $g \colon R \to R$ for which
   \begin{equation}
   \label{eq:endbar}
   g(r)^p ~\equiv~ g(r^p) \quad \pmod{pR}
   \end{equation}
for all primes $p$ and elements $r \in R$.  This last condition is clearly equivalent to
   \begin{equation}
   \label{eq:endbar p}
   g(x^i)^p ~\equiv~ g(x^{ip}) \quad \pmod{pR}
   \end{equation}
for $i = 0, 1$ and all primes $p$.  Consider first the case $i = 0$.  Writing $g(1) = a + bx$, the condition then says
   \[
   \begin{split}
   a + bx 
   &~\equiv~ (a + bx)^p \pmod{p} \\
   &~\equiv~ a + bx^p \pmod{p} \\
   &~=~ a,
   \end{split}
   \]
since $x^2 = 0$ in $R$.  In other words, $b = 0$ and $g(1) = a \in \bZ$.  Now write $g(x) = c + dx$.  Since $x^p = 0$ for any prime $p$ and since $g(0) = 0$, the condition above for $i = 1$ says
   \[
   \begin{split}
   g(x^p) &~=~ 0 \\
   &~\equiv~ (c + dx)^p \pmod{p} \\
   &~\equiv~ c + dx^p \pmod{p} \\
   &~=~ c.
   \end{split}
   \]
Since this is true for all primes $p$, we infer that $c = 0$ and $g(x) = dx$.

Summarizing this discussion, we have shown that $\Endbar(R)$ consists of precisely the $\bZ$-linear self-maps $g$ of $R$ for which $g(1) \in \bZ$ and $g(x) \in \bZ x$.  It is clear that any such map commutes with $\psi^k$ for all $k$.  Therefore, $\Hl^0(R) = \Endbar(R)$ is isomorphic to the product ring $\bZ \oplus \bZ$, as stated in the statement of the Theorem.

\emph{Step 2}. $\Hl^1(R)$.

Observe that, since each $g \in \cF^0(R) = \Endbar(R)$ commutes with all the $\psi^k$ in $R$, the image of the differential
   \[
   d^0 ~=~ \lbrack \psi^*,\, - \rbrack \colon \cF^0(R) ~\to~ \cF^1(R)
   \]
is trivial.  Therefore, $\Hl^1(R) = \Der(R)$, and so it suffices to compute $\Derbar(R)$.

To do this, let $(f(p))_p$ be a sequence of $\bZ$-linear self maps of $R$ indexed by the primes with $f(p)(R) \subseteq pR$.  Thus, using the $\bZ$-basis $\lbrace 1, x \rbrace$ of $\bZ \lbrack x \rbrack/(x^2)$, we can write each $f(p)$ as a $2$-by-$2$ matrix with entries in $p\bZ$, say, $f(p) = p(a(p)_{ij})$.  If $\psi^p(x) = b_px$ in $R$, then we can similarly represent $\psi^p$ as a $2$-by-$2$ diagonal matrix with entries $1$ and $b_p$ along the diagonal.  In this context, we have
   \[
   \begin{split}
   \psi^pf(q) ~+~& f(p)\psi^q \\
   &~=~
   q\begin{bmatrix} 1 & \\ & b_p \end{bmatrix}
   \begin{bmatrix} a(q)_{11} & a(q)_{12} \\ a(q)_{21} & a(q)_{22} \end{bmatrix} 
   ~+~
   p \begin{bmatrix} a(p)_{11} & a(p)_{12} \\ a(p)_{21} & a(p)_{22} \end{bmatrix}
   \begin{bmatrix} 1 & \\ & b_q \end{bmatrix} \\
   &~=~ \begin{bmatrix} 
        pa(p)_{11} + qa(q)_{11}    & pb_qa(p)_{12} + qa(q)_{12} \\
        pa(p)_{21} + qb_pa(q)_{21} & pb_qa(p)_{22} + qb_pa(q)_{22} 
        \end{bmatrix}
   \end{split}
   \]
As usual, the empty entries denote $0$.  One obtains the matrix representation of $(\psi^qf(p) + f(q)\psi^p)$ by interchanging $p$ and $q$ in the above matrix.  Using this, the condition \eqref{eq:f(pq)} can now be seen to be equivalent to the equalities
   \[
   pa(p)_{ij}(b_q - 1) ~=~ qa(q)_{ij}(b_p - 1) 
   \]
for all primes $p$ and $q$ and $(i, j) = (1, 2)$ and $(2, 1)$.  Since each $b_p$ is divisible by $p$, it follows that $a(q)_{ij}$ is divisible by $p$ for all $p \not= q$, i.e.\ $a(q)_{ij} = 0$.  It follows that $a(p)_{ij} = 0$ for all primes $p$ and  $(i, j) = (1, 2)$ and $(2, 1)$.  In other words, each $f(p)$ is a diagonal matrix with entries in $p \bZ$.  Therefore, we have the group isomorphisms
   \[
   \Hl^1(R) ~\cong~ \Derbar(R) ~\cong~ \prod_{p \in \cP}\, (p \bZ \oplus p \bZ).
   \]

\emph{Step 3}: Ring structure.

To finish the proof, we only need to observe that if $(a, b) \in \Hl^0(R)$ and $\lbrace (x_p, y_p) \rbrace \in \Hl^1(R)$, then, by Corollary \ref{cor1:H1}, we have
   \[
   \begin{split}
   (a,\, b) \circ \lbrace (x_p,\, y_p) \rbrace
    &~=~ \left\lbrace \begin{bmatrix} a & \\ & b \end{bmatrix}
        \begin{bmatrix} x_p & \\ & y_p \end{bmatrix} \right\rbrace
    ~=~ \left\lbrace\begin{bmatrix} ax_p & \\  & by_p \end{bmatrix}\right\rbrace \\
    &~=~ \lbrace (ax_p,\, by_p) \rbrace
    ~=~ \lbrace (x_p,\, y_p) \rbrace \circ (a,\, b),
   \end{split}
   \]
as desired.
\end{proof}


\section{$\Hl^0$ of filtered $\lambda$-ring structures on $\bZ \lbrack x \rbrack/(x^3)$}
\label{sec:H0 n=3}

In this section, we will compute the algebra $\Hl^0(R)$ for each of the uncountably many isomorphism classes of filtered $\lambda$-ring structures on $\bZ \lbrack x \rbrack/(x^3)$, with $x$ in some fixed positive filtration.  (See \cite{yau4} for a proof of this uncountability statement.)


\subsection{Filtered $\lambda$-ring structures on $\bZ \lbrack x \rbrack/(x^3)$}
\label{subsec1:n=3}

Let us begin by recalling the classification of filtered $\lambda$-ring structures on $\bZ \lbrack x \rbrack/(x^3)$.  See \cite{yau4} for more details.

In what follows, we will describe a filtered $\lambda$-ring structure on $\bZ \lbrack x \rbrack/(x^3)$ in terms of its Adams operations $\psi^p$ for $p$ primes.  This is sufficient to determine the filtered $\lambda$-ring structure by a result of Wilkerson \cite{wil}, as discussed in section \ref{subsec:adams}.

Let $R$ be a filtered $\lambda$-ring structure on $\bZ \lbrack x \rbrack/(x^3)$ with Adams operations
   \begin{equation}
   \label{eq:psi}
   \psi^p_R(x) ~=~ \psi^p(x) ~=~ b_px + c_px^2.
   \end{equation}
It is shown in \cite{yau4} that if $b_2 = 0$, then 
   \begin{itemize}
   \item $b_p = 0$ for all primes $p$,
   \item $c_2$ is an odd integer, and 
   \item $c_p \equiv 0 \pmod{p}$ for all odd primes $p$.  
   \end{itemize}
In this case, we write $S((c_p))$ for $R$, since its filtered $\lambda$-ring structure is completely determined by the $c_p$ for $p$ primes.  Conversely, any such sequence $(c_p)$ gives rise to a filtered $\lambda$-ring structure on $\bZ \lbrack x \rbrack/(x^3)$.  Two such filtered $\lambda$-ring structures, $S((c_p))$ and $S((c_p^\prime))$, are isomorphic if and only if $(c_p) = \pm(c_p^\prime)$.  In particular, there are uncountably many isomorphism classes of filtered $\lambda$-ring structures on $\bZ \lbrack x \rbrack/(x^3)$ of the form $S((c_p))$.

On the other hand, if $b_2 \not = 0$, then there exists an odd integer $h$ such that
   \begin{equation}
   \label{eq:cp}
   c_p ~=~ h\frac{b_p(b_p - 1)}{G}
   \end{equation}
for all primes $p$, where 
   \[
   G ~=~ \gcd \lbrace b_q(b_q - 1) \colon \text{ all primes }q \rbrace.
   \]
Moreover, the following conditions have to hold: 
   \begin{itemize}
   \item $b_p \equiv 0 \pmod{p}$ for all primes $p$.
   \item $b_p(b_p - 1) \equiv 0 \pmod{2^{\nu_2(b_2)}}$, where $\nu_2$ denotes the $2$-adic evaluation of an integer (i.e.\ the exponent of the prime factor $2$ in the integer).
   \item $h$ is in the range $1 \leq h \leq G/2$.
   \item Suppose that there exists an odd prime $p$ for which $b_p \not= 0$ and
   \[
   \nu_p(b_p) ~=~ \min \lbrace \nu_p(b_q(b_q - 1)) \colon \, b_q \not= 0 \rbrace.
   \]
(There are at most finitely many such primes, since each such $p$ divides $b_2(b_2 - 1) \not= 0$.)  Then any such prime $p$ divides $h$.
   \end{itemize}
In this case, we denote $R$ by $S((b_p), h)$, since the entire filtered $\lambda$-ring structure is determined by the $b_p$ for $p \in \cP$ and $h$.  Conversely, given integers $b_p$ ($p \in \cP$) and $h$ satisfying the above properties, there is a filtered $\lambda$-ring structure on $\bZ \lbrack x \rbrack/(x^3)$ whose $\psi^p$ is given by \eqref{eq:psi}, in which $c_p$ is determined by $b_p$ and $h$ via \eqref{eq:cp}.   Two such filtered $\lambda$-rings, $S((b_p), h)$ and $S((b_p^\prime), h^\prime)$, are isomorphic if and only if $b_p = b_p^\prime$ for all primes $p$ and $h = h^\prime$.  In particular, there are uncountably many isomorphism classes of filtered $\lambda$-ring structures on $\bZ \lbrack x \rbrack/(x^3)$ of the form $S((b_p), h)$.

For example, let $\mathbf{F}$ denote the complex numbers $\mathbf{C}$, the quaternions $\mathbf{H}$, or the Cayley octonions $\mathbf{O}$, and let $\mathbf{FP}^2$ be the corresponding projective $2$-space.  Then the  $K$-theory filtered $\lambda$-ring of $\mathbf{FP}^2$ is
   \begin{equation}
   \label{eq:proj}
   K(\mathbf{FP}^2) ~=~ \begin{cases} 
   S((p),\, 1) \text{ with } G = 2(2 - 1) = 2 & \text{ if } \mathbf{F} = \mathbf{C} \\
   S((p^2), \, 1) \text{ with } G = 2^2(2^2 - 1) = 12 & \text{ if } \mathbf{F} = \mathbf{H} \\
   S((p^4), \, 1) \text{ with } G = 2^4(2^4 - 1) = 240 & \text{ if } \mathbf{F} = \mathbf{O}.\end{cases}
   \end{equation}


\subsection{The algebras $\Hl^0(R)$}
\label{subsec2:n=3}

Let $R$ be any filtered $\lambda$-ring structure on $\bZ \lbrack x \rbrack/(x^3)$.  Using $\lbrace 1, x, x^2 \rbrace$ as a $\bZ$-basis for $\bZ \lbrack x \rbrack/(x^3)$, we can write each $\bZ$-linear self map of $R$ as an element of $M(3, \bZ)$, the algebra of $3$-by-$3$ matrices with integer entries.  We will continue to omit entries that are $0$ in a matrix.

Here is the main result of this section.  Recall that $\cP$ denotes the set of all primes.

\begin{thm}
\label{thm:H0 n=3}
Let $R$ be a filtered $\lambda$-ring structure on $\bZ \lbrack x \rbrack/(x^3)$. 
   \begin{enumerate}
   \item If $R = S((c_p))$ for some $c_p$, $p \in \cP$, then 
   \[
   \Hl^0(R) ~=~ \left\lbrace 
                \begin{bmatrix} a &  & \\
                                 & j & \\
                                 & k & j
                \end{bmatrix} ~\colon~~ a, \, j,\, k \in \bZ \right\rbrace
   \]
as a subalgebra of $M(3, \bZ)$.
   \item If $R = S((b_p), h)$ for some $b_p$, $p \in \cP$, and $h$, then 
   \[
   \Hl^0(R) ~=~ \left\lbrace
                \begin{bmatrix}
                a &  & \\
                 & j &  \\
                 & k & t \end{bmatrix} ~\colon~~ h(t - j) = kG \right\rbrace
   \]
as a subalgebra of $M(3, \bZ)$, where $G = \gcd(b_q(b_q - 1))_{q \in \cP}$.
   \end{enumerate}
In each case, $\Hl^0(R)$ is a commutative algebra and has rank $3$ over $\bZ$ as an additive group.
\end{thm}

To prove this Theorem, we will first compute $\cF^0(R) = \Endbar(R)$, still using the $\bZ$-basis $\lbrace 1, x, x^2 \rbrace$.

\begin{lemma}
\label{lem1:H0 n=3}
Let $R$ be any filtered $\lambda$-ring structure on $\bZ \lbrack x \rbrack/(x^3)$ and let $g$ be a $\bZ$-linear self map of $R$.  Then $g \in \Endbar(R) = \cF^0(R)$ if and only if 
   \[
   g ~=~ \begin{bmatrix}
         a &  & \\
         & j & s \\
         & k & t \end{bmatrix}
   \]
with
   \[
   s ~\equiv~ 0 ~\equiv~ t - j \pmod{2}.
   \]
\end{lemma}

\begin{proof}
By definition, the $\bZ$-linear map $g$ lies in $\Endbar(R)$ if and only if \eqref{eq:endbar} holds for all primes $p$ and elements $r \in R$.  This condition only needs to be checked for $r = 1$, $x$, and $x^2$.  If $g(1) = a + bx + cx^2$, then \eqref{eq:endbar}, when applied to $r = 1$, becomes
   \[
   \begin{split}
   g(1^p) &~=~ g(1) ~=~ a + bx + cx^2 \\
          &~\equiv~ (a + bx + cx^2)^p \pmod{p} \\
          &~\equiv~ a + bx^p + cx^{2p} \pmod{p} \\
          &~=~ a \pmod{p} \text{ if } p > 2.
   \end{split}
   \]
It follows that $b = c = 0$.  The case $p = 2$ gives no additional information.

Write $g(x) = i + jx + kx^2$.  Then \eqref{eq:endbar}, when applied to $r = x$, becomes
   \[
   \begin{split}
   g(x^p) &~=~ 0 \text{ if } p > 2 \\
          &~\equiv~ (i + jx + kx^2)^p \pmod{p} \\
          &~\equiv~ i \pmod{p}.
   \end{split}
   \]
Since this is true for all odd primes $p$, we have that $i = 0$.  The case $p = 2$ has not been used yet; we will come back to this below.

Write $g(x^2) = r + sx + tx^2$.  Then \eqref{eq:endbar} becomes
   \[
   \begin{split}
   0 &~=~ g(x^{2p}) \\
     &~\equiv~ (r + sx + tx^2)^p \pmod{p} \\
     &~\equiv~ r + sx^p \pmod{p}.
   \end{split}
   \]
Therefore, $r = 0$ and $s \equiv 0 \pmod{2}$.  Moreover, the condition
   \[
   g(x)^2 ~\equiv~ g(x^2) \pmod{2}
   \]
is equivalent to
   \[
   jx^2 ~\equiv~ tx^2 \pmod{2},
   \]
i.e.\ $j \equiv t \pmod{2}$.

This finishes the proof of the Lemma.
\end{proof}

\begin{proof}[Proof of Theorem \ref{thm:H0 n=3}]
It is immediate from the definition that $\Hl^0(R)$ is the subalgebra of $\Endbar(R)$ consisting of those $\bZ$-linear self maps $g$ of $R$ for which $g\psi^p = \psi^p g$ for all primes $p$.  Write $\psi^p(x) = b_px + c_px^2$ and let $g \in \Endbar(R)$ be as in the statement of Lemma \ref{lem1:H0 n=3}.  Then the equation $g \psi^p ~=~ \psi^p g$ 
can be written in matrix form as
   \begin{equation}
   \label{eq1:H0 n=3}
   g \psi^p ~=~
   \begin{bmatrix}
   a & & \\
   & jb_p + sc_p & sb_p^2 \\
   & kb_p + tc_p & tb_p^2
   \end{bmatrix}
   ~=~
   \begin{bmatrix}
   a & & \\
   & jb_p & sb_p \\
   & jc_p + kb_p^2 & sc_p + tb_p^2
   \end{bmatrix}
    ~=~ \psi^p g
   .
   \end{equation}
In the case that $R = S((c_p))$, $b_p = 0$ for all $p$.  Therefore, \eqref{eq1:H0 n=3} is equivalent to the conditions
   \[
   sc_p ~=~ 0 ~=~ (t - j)c_p.
   \]
Since $c_2$ is an odd integer, it is non-zero in any case.  It follows that
   \[
   s ~=~ 0 ~=~ t - j,
   \]
and the first part of Theorem \ref{thm:H0 n=3} is proved.

Now suppose that $R = S((b_p), h)$ for some $b_p$ and $h$.  By comparing the $(2,\, 3)$ entries in \eqref{eq1:H0 n=3}, using the fact that $b_2 \not= 0$ is an even integer, one infers that $s = 0$.  The only condition left in \eqref{eq1:H0 n=3} now is
   \[
   kb_p + tc_p ~=~ jc_p + kb_p^2,
   \]
or equivalently,
   \begin{equation}
   \label{eq2:H0 n=3}
   h\frac{b_p(b_p - 1)}{G}(t - j) ~=~ kb_p(b_p - 1).
   \end{equation}
This condition is trivially true if $b_p = 0$.  If $b_p \not= 0$, then, since $b_p \not=  1$ in any case, this condition is equivalent to
   \begin{equation}
   \label{eq:hkG}
   h(t - j) ~=~ kG.
   \end{equation}
Conversely, it is easy to see that the conditions, $s = 0$ and \eqref{eq:hkG}, imply \eqref{eq1:H0 n=3} for all primes $p$.

This proves the second part of Theorem \ref{thm:H0 n=3}.

It remains to establish the last statement in the Theorem.  When $R = S((c_p))$, it is easy to see that $\Hl^0(R)$ is free of rank $3$ over $\bZ$ as a group.  Indeed, the following elements form a $\bZ$-basis for $\Hl^0(R)$:  
   \[
   \begin{bmatrix}
   1 & & \\
     & & \\
     & & 
   \end{bmatrix},\,
   \begin{bmatrix}
    & & \\
    & 1 & \\
    & & 1 
   \end{bmatrix},\,
   \begin{bmatrix}
    & & \\
    & & \\
    & 1 & 
   \end{bmatrix}.
   \]
Moreover, any two elements in $\Hl^0(S((c_p)))$ commute, since
   \[
   \begin{bmatrix} a &  & \\
                     & j & \\
                     & k & j
   \end{bmatrix}
   \begin{bmatrix} a^\prime &  & \\
                     & j^\prime & \\
                     & k^\prime & j^\prime
   \end{bmatrix}
   ~=~
   \begin{bmatrix} aa^\prime &  & \\
                     & jj^\prime & \\
                     & kj^\prime + jk^\prime & jj^\prime
   \end{bmatrix}.
   \]
This last matrix remains the same if $a$ (respectively $j$ and $k$) and $a^\prime$ (respectively $j^\prime$ and $k^\prime$) are interchanged.  This shows that the algebra $\Hl^0(S((c_p)))$ is commutative.

For the second case, $R = S((b_p), h)$, we have
   \[   
   \begin{bmatrix}  
   a &   & \\
   & j & \\
   & k & t \end{bmatrix} 
   \begin{bmatrix}  
   a^\prime &   & \\
   & j^\prime & \\
   & k^\prime & t^\prime \end{bmatrix} 
   ~=~
   \begin{bmatrix}
   aa^\prime & & \\
   & jj^\prime & \\
   & kj^\prime + tk^\prime & tt^\prime
   \end{bmatrix}
   \]
To see that this is equal to the product with the reserve order, observe that we can write $kk^\prime G/h$ in two ways:
   \[
   \frac{kk^\prime G}{h} ~=~ k^\prime(t - j) ~=~ k(t^\prime - j^\prime).
   \]
Therefore, we have that
   \[
   kj^\prime + tk^\prime ~=~ k^\prime j + t^\prime k,
   \]
which shows that $\Hl^0(R)$ is a commutative algebra.  To see that $\Hl^0(R)$ is free of rank $3$ over $\bZ$, observe that 
   \begin{equation}
   \label{eq:jkt}
   \lbrace (j, k, t) \in \bZ^{\times 3} \colon h(t - j) = kG \rbrace
   \end{equation}
is the kernel of the surjective $\bZ$-linear map
   \[
   \varphi \colon \bZ^{\times 3} ~\to~ d\bZ ~\cong~ \bZ,
   \]
where $d = gcd(h, G)$ and 
   \[
   \varphi((j, k, t)) ~=~ h(t - j) - kG.
   \]
So the group in \eqref{eq:jkt} is free of rank $2$.  It follows that $\Hl^0(R)$ is free of rank $3$, as the $(1, 1)$ entries of elements in it are arbitrary.

This finishes the proof of Theorem \ref{thm:H0 n=3}.
\end{proof}


\section{$\Hl^1$ of the $64$ $S((p^r), h)$}
\label{sec:H1 n=3}

The main purpose of this section is to compute the groups $\Hl^1(R)$ and to determine the commutativity of the graded algebras $\Hl^{\leq 1}(R)$ for the following $64$ filtered $\lambda$-ring structures $R$ over $\bZ \lbrack x \rbrack/(x^3)$:
   \begin{equation}
   \label{eq:S}
   S((p^r),\, h) ~=~ \left \lbrace \psi^p(x) ~=~ p^rx ~+~ h\frac{p^r(p^r - 1)}{2^r(2^r - 1)}x^2 \right\rbrace.
   \end{equation}
Here $r \in \lbrace 1, 2, 4 \rbrace$ and 
   \[
   h ~\in~ \begin{cases} 
           \lbrace 1 \rbrace & \text{ if } r = 1, \\
           \lbrace 1, \, 3,\, 5 \rbrace & \text{ if } r = 2, \\
           \lbrace 1,\, 3,\, \ldots ,\, 119 \rbrace & \text{ if } r = 4.
           \end{cases}
   \]
In the notation of section \ref{subsec1:n=3}, these are the only filtered $\lambda$-ring structures on $\bZ \lbrack x \rbrack/(x^3)$ (up to isomorphism) of the form $S((p^r),\, h)$ for $r = 1$, $2$, and $4$.

Before we proceed, we should explain the significance of these $64$ filtered $\lambda$-rings.  The following statement is shown in \cite{yau4}:  If $X$ is a torsionfree topological space (i.e.\ its integral cohomology is $\bZ$-torsionfree) whose unitary $K$-theory $K(X)$, as a filtered ring, is the ring $\bZ \lbrack x \rbrack/(x^3)$, then $K(X)$ is isomorphic as a filtered $\lambda$-ring to one of the $64$ $S((p^r), h)$ in \eqref{eq:S}.  
In other words, among the uncountably many isomorphism classes of filtered $\lambda$-ring structures on $\bZ \lbrack x \rbrack/(x^3)$, only these $64$ isomorphism classes can possibly be topologically realized by torsionfree spaces.

Given one of these $64$ $\lambda$-rings $S((p^r), h)$, define 
   \begin{equation}
   \label{eq:D}
   D ~=~ D(r, h) ~\buildrel \text{def} \over =~ \gcd(h, 2^r(2^r - 1)).
   \end{equation}
For example, $D = 1$ if $h = 1$.  These are the three cases for the $K$-theory of the projective $2$-spaces $\mathbf{FP}^2$ with $\mathbf{F} = \mathbf{C}$ $(r = 1)$, $\mathbf{H}$ $(r = 2)$, and $\mathbf{O}$ $(r = 4)$ (see \eqref{eq:proj}).  Exactly $35$ of these $D(r, h)$ are equal to $1$, namely, $D(1, 1)$, $D(2, 1)$, $D(2, 5)$, and $32$ of the $D(4, h)$.  Also, let $\bZ^{\bN}$ denote a countably infinite product of copies of the additive group of integers $\bZ$.

With these notations, we can now state the main result of this section.

\begin{thm}
\label{thm:H1 n=3}
Let $R = S((p^r), h)$ be any one of the $64$ filtered $\lambda$-rings in \eqref{eq:S}.  Then there is a group isomorphism
   \[
   \Hl^1(R) ~\cong~ \frac{\bZ}{h\bZ} ~\times~ \frac{\bZ}{D\bZ} ~\times~ \bZ^{\bN}.
   \]
The graded algebra $\Hl^{\leq 1}(R)$ is commutative if and only if $D = 1$.
\end{thm}

We will need an explicit description of $\Derbar(R)$.  Using the $\bZ$-basis $\lbrace 1, x, x^2 \rbrace$, if $f$ is a $\bZ$-linear self map of $R$ satisfying $f(R) \subseteq pR$ for some $p$, then it can be represented by a $3$-by-$3$ matrix $(pa_{ij})$ with entries in $p \bZ$.

\begin{lemma}
\label{lem1:H1 n=3}
Let 
   \[
   R ~=~ S((b_p),\, h) ~=~ \lbrace \psi^p(x) = b_px + c_px^2 \rbrace
   \]
be any filtered $\lambda$-ring structure on $\bZ \lbrack x \rbrack/(x^3)$ with $b_2 \not= 0$.  Let $\lbrace f(p) \rbrace_{p \in \cP} = \lbrace (pa(p)_{ij}) \rbrace_{p \in \cP}$ be a sequence of $\bZ$-linear self maps of $R$ with $f(p)(R) \subseteq pR$ for each $p$.  Then $\lbrace f(p) \rbrace_{p \in \cP}$ is an element of $\Derbar(R)$ if and only if the following three conditions hold for all primes $p$ and $q$:
   \begin{equation}
   \label{eq:0}
   a(p)_{12} = a(p)_{13} = a(p)_{21} = a(p)_{31} = 0,
   \end{equation}
   \begin{equation}
   \label{eq:23}
   pb_q(b_q - 1)a(p)_{23} ~=~ qb_p(b_p - 1)a(q)_{23},
   \end{equation}
and
   \begin{equation}
   \begin{split}
   \label{eq:32}
   & \quad qc_p\left(a(q)_{22} ~-~ a(q)_{33}\right) ~+~ qb_p(b_p ~-~ 1)a(q)_{32} \\
    &~=~ pc_q\left(a(p)_{22} ~-~ a(p)_{33}\right) ~+~ pb_q(b_q ~-~ 1)a(p)_{32}.
   \end{split}
   \end{equation}
\end{lemma}

\begin{proof}
A little bit of matrix computation shows that for any primes $p$ and $q$, one has
\begin{small}
   \[
   \begin{split}
   & \psi^p f(q) ~+~ f(p)\psi^q ~=~ \\
   & 
   \begin{bmatrix}
   pa(p)_{11} + qa(q)_{11} & qa(q)_{12} + pb_qa(p)_{12} + pc_qa(p)_{13} & qa(q)_{13} + pb_q^2a(p)_{13} \\
   & & \\
   qb_pa(q)_{21} + pa(p)_{21} & qb_pa(q)_{22} + pb_qa(p)_{22} + pc_qa(p)_{23} & qb_pa(q)_{23} + pb_q^2a(p)_{23} \\
   & & \\
   \begin{split} 
   & qc_pa(q)_{21} + \\ 
   & qb_p^2a(q)_{31} + pa(p)_{31}
   \end{split} 
   & 
   \begin{split}
   & qc_pa(q)_{22} + qb_p^2a(q)_{32} \\
   &+ pb_qa(p)_{32} + pc_qa(p)_{33} 
   \end{split}
   & 
   \begin{split}
   & qc_pa(q)_{23} + \\
   & qb_p^2a(q)_{33} + pb_q^2a(p)_{33}
   \end{split}
   \end{bmatrix}
   \end{split}
   \]
\end{small}

The matrix for the map $(\psi^q f(p) + f(q) \psi^p)$ is obtained from this matrix by interchanging $p$ and $q$.  The sequence $\lbrace f(p) \rbrace$ is an element of $\Derbar(R)$ if and only if \eqref{eq:f(pq)} holds, which is equivalent to saying that the $3$-by-$3$ matrices represented by the two sides are equal.  We will now think of \eqref{eq:f(pq)} as a matrix equation.  It is clear that the $(1,\, 1)$ entries give no information, and the $a(p)_{11}$ can be any integers.

The equality of the $(1, 3)$ entries in \eqref{eq:f(pq)} is equivalent to
   \[
   qa(q)_{13} ~+~ pb_q^2a(p)_{13} ~=~ pa(p)_{13} ~+~ qb_p^2a(q)_{13},
   \]
which can be rewritten as 
   \[
   p(b_q^2 - 1)a(p)_{13} ~=~ q(b_p^2 - 1)a(q)_{13}.
   \]
Since each $b_p$ is a multiple of $p$, this equation forces $a(p)_{13} = 0$ for all $p$.  Applying this to the $(1,\, 2)$ entries in \eqref{eq:f(pq)}, one obtains
   \[
   q(b_p - 1)a(q)_{12} ~=~ p(b_q - 1)a(p)_{12}.
   \]
Just as above, this forces $a(p)_{12} = 0$ for all $p$.  A similar argument applies to the $(2, 1)$ and the $(3,\, 1)$ entries in \eqref{eq:f(pq)} and gives rise to $a(p)_{21} = a(p)_{31} = 0$ for all primes $p$.  In other words, the condition imposed by the first rows and the first columns in \eqref{eq:f(pq)} is exactly \eqref{eq:0}.

From the matrix above, it is immediate that the equality of the $(2, 3)$ (respectively $(3, 2)$) entries in \eqref{eq:f(pq)} is exactly \eqref{eq:23} (respectively \eqref{eq:32}).

It remains to show that the $(2, 2)$ and the $(3, 3)$ entries in \eqref{eq:f(pq)} give rise to redundant conditions.  In each case, the condition is
   \[
   pc_qa(p)_{23} ~=~ qc_pa(q)_{23}.
   \]
This can be obtained from \eqref{eq:23} by multiplying both sides of that equation by $h/G$ (see \eqref{eq:cp}), as claimed.  (It is here that we are using the condition $b_2 \not= 0$.)

This finishes the proof of the Lemma.
\end{proof}

We also need an explicit description of $\Innbar(R)$.

\begin{lemma}
\label{lem2:H1 n=3}
Let $R$ be any filtered $\lambda$-ring structure on $\bZ \lbrack x \rbrack/(x^3)$ with $\psi^p(x) = b_p x + c_px^2$.  Then, using the $\bZ$-basis $\lbrace 1,  x, x^2 \rbrace$,  $\Innbar(R)$ consists of the following sequences (indexed by the primes) of matrices,
   \begin{equation}
   \label{eq:inn}
   \left \lbrace 
   \begin{bmatrix} 0 & 0 & 0\\        
   0 & -sc_p & -sb_p(b_p - 1) \\
   0 & kb_p(b_p - 1) + (j - t)c_p & sc_p \end{bmatrix}\right\rbrace_{p \in \cP},
   \end{equation}
such that $s \equiv 0 \equiv t - j \pmod{2}$.
\end{lemma}

Notice that any such sequence of matrices is completely determined by the four parameters $j$, $k$, $s$, and $t$.

\begin{proof}
This is immediate from Lemma \ref{lem1:H0 n=3} and \eqref{eq1:H0 n=3}, noting that the matrix displayed above is exactly $\lbrack \psi^p, g \rbrack$ with $g$ as in the statement of Lemma \ref{lem1:H0 n=3}.
\end{proof}


\subsection{Proof of Theorem \ref{thm:H1 n=3} when $r = 1$}
\label{subsec:r=1}

Using Lemma \ref{lem1:H1 n=3}, we infer that $\Derbar(R)$ $=$ $\Derbar(S((p), 1)$ consists of the sequences $\lbrace f(p) = p(a_{ij}) \rbrace_{p \in \cP}$ of $\bZ$-linear self maps of $R$ satisfying \eqref{eq:0},
   \begin{equation}
   \label{eq:23 r=1}
   a(p)_{23} ~=~ (p - 1)a(2)_{23},
   \end{equation}
and 
   \begin{equation}
   \label{eq:32 r=1}
   \begin{split}
   a(p)_{22} &~-~ a(p)_{33} ~+~ 2a(p)_{32} \\
   &~=~ (p - 1)\bigl(a(2)_{22} ~-~ a(2)_{33} ~+~ 2a(2)_{32}\bigr)
   \end{split}
   \end{equation}
for all primes $p$.  In other words, to obtain an element in $\Derbar(R)$, one  chooses five arbitrary elements in $2 \bZ$ for the $(1, 1)$, $(2, 2)$, $(2, 3)$, $(3, 2)$, and $(3, 3)$ entries in $f(2)$.  Then for each odd prime $p$, one chooses three arbitrary elements in $p \bZ$ for the $(1, 1)$, $(2, 2)$, and $(3, 2)$ entries in $f(p)$.  The $(2, 3)$ and $(3, 3)$ entries are then determined by \eqref{eq:23 r=1} and \eqref{eq:32 r=1}, respectively.  Therefore, we can make the identification
   \begin{equation}
   \label{eq:der r=1}
   \Derbar(R) ~=~ (2 \bZ)^{\times 5} ~\times \prod_{p > 2}\, (p\bZ)^{\times  3}.
   \end{equation}

As for $\Innbar(R)$, note that the element displayed in \eqref{eq:inn} in Lemma \ref{lem2:H1 n=3} now becomes
   \[
   \left \lbrace
   \binom{p}{2}
   \begin{bmatrix}
   0 & 0 & 0 \\
   0 & -s & -2s \\
   0 & j - t + 2k & s
   \end{bmatrix}\right \rbrace_{p \in \cP}.
   \]
In particular, the $p = 2$ component of this element  is the matrix
   \[
   \begin{bmatrix}
   & & \\
   & -s & -2s \\
   & j - t + 2k & s
   \end{bmatrix}.
   \] 
As these elements vary through $\Innbar(R)$, $s$ and $(j - t)$ run through all even integers.  Therefore, using the identification \eqref{eq:der r=1}, we conclude that
   \[
   \Hl^1(R) 
   ~\cong~ \frac{\Derbar(R)}{\Innbar(R)} 
   ~\cong~ (2\bZ)^{\times 3} \times \prod_{p > 2}\, (p\bZ)^{\times  3} 
   ~\cong~ \bZ^{\bN}
   \]
as groups.  This is exactly the assertion of Theorem \ref{thm:H1 n=3} in the case $r = 1$, since $h = D = 1$ here.

The commutativity assertion is dealt with below, after the $r = 4$ case.


\subsection{Proof of Theorem \ref{thm:H1 n=3} when $r = 2$}
\label{subsec:r=2}

In this case, \eqref{eq:23} and \eqref{eq:32} become, respectively,
   \begin{equation}
   \label{eq:23 r=2}
   a(p)_{23} ~=~ \frac{p(p^2 - 1)}{6} a(2)_{23}
   \end{equation}
and
   \begin{equation}
   \label{eq:32 r=2}
   \begin{split}
   h(a(p)_{22} &~-~ a(p)_{33}) ~+~ 12 a(p)_{32} \\
   &~=~ \frac{p(p^2 - 1)}{6}(h(a(2)_{22} ~-~ a(2)_{33}) ~+~ 12 a(2)_{32}),
   \end{split}
   \end{equation}
where $R = S((p^2), h)$ with $h \in \lbrace 1,  3,  5 \rbrace$.  If $h = 1$ or $3$, then the same argument as in the $r = 1$ case allows us to once again make the identification \eqref{eq:der r=1} for $\Derbar(R)$.

If $h = 5$, then, since $5$ and $12$ are relatively prime, the left-hand side of \eqref{eq:32 r=2} can still attain any integer by choosing $a(p)_{22}$, $a(p)_{33}$, and $a(p)_{32}$ appropriately.   In particular, to choose an element $\lbrace f(p) \rbrace$ of $\Derbar(R)$, one can choose five arbitrary elements in $2\bZ$ for the $(1, 1)$, $(2, 2)$, $(2, 3)$, $(3, 2)$, and $(3, 3)$ entries in $f(2)$.  Furthermore, for each odd prime $p$, one chooses an element in $p \bZ$ for the $(1, 1)$ entry in $f(p)$ and two more arbitrary integers to form $f(p)$.  Indeed, denoting the value in the right-hand side of \eqref{eq:32 r=2} by $N_p$, one has
   \begin{equation}
   \label{eq:freedom}
   \begin{split}
   a(p)_{32} &~=~ -2N_p ~+~ 5r_p \\
   a(p)_{22} ~-~ a(p)_{33} &~=~ 5N_p ~-~ 12r_p
   \end{split}
   \end{equation}
for some integer $r_p$.  Therefore, there are two degrees of freedom left, namely, $r_p$ and, say, $a(p)_{22}$.  This allows us to make the identification
   \begin{equation}
   \label{eq:der r=2 h=5}
   \begin{split}
   \Derbar(S((p^2),\, 5) 
   &~=~ (2\bZ)^{\times 5} ~\times~
   \prod_{p > 2}\,  (p\bZ \times \bZ \times \bZ) \\
   &~\cong~  (2\bZ)^{\times 5} ~\times~ \bZ^{\bN}.
   \end{split}
   \end{equation}

As for $\Innbar(R)$, note once again that each element in $\Innbar(R)$ is still determined by the four parameters $s$, $j$, $k$, and $t$ with $s$ and $(j - t)$ even.  The $p = 2$ component of a typical element  in  $\Innbar(R)$ (see \eqref{eq:inn}) now takes the form
   \[
   \begin{bmatrix}
   & & \\
   & -sh & -12 s \\
   & h(j - t) + 12k & sh
   \end{bmatrix}.
   \]
To finish the proof,  observe that
   \[
   \lbrace sh \colon s \in 2\bZ \rbrace ~=~ 2h\bZ
   \]
and
   \[
   \lbrace  h(j - t) + 12k \colon (j - t) \in 2 \bZ,\, k \in \bZ \rbrace 
   ~=~ \begin{cases} 2\bZ & \text{if } h = 1,\, 5, \\
                                   6\bZ & \text{if } h = 3 \end{cases}.
   \]  
 Summarizing this discussion, we have shown that
   \[
   \Hl^1(S((p^2),\, h)) 
   ~\cong~ 
   \begin{cases}
   (2\bZ)^{\times 3} \times \bZ^{\bN} &\text{if } h = 1\\
   (2\bZ/6\bZ)^{\times 2} \times (2\bZ)^{\times 3} \times \bZ^{\bN} & \text{if } h = 3 \\
   (2\bZ/10 \bZ) \times (2\bZ)^{\times 3} \times \bZ^{\bN} & \text{if } h = 5\end{cases}
   \]
This is exactly the claim in Theorem \ref{thm:H1 n=3} when $r = 2$, since $D = 1$ (respectively $3$) when $h = 1$ or $5$ (respectively $3$).

The commutativity assertion is dealt with below, after the $r = 4$ case.


\subsection{Proof of Theorem \ref{thm:H1 n=3} when $r = 4$}
\label{subsec:r=4}

When $r = 4$,  \eqref{eq:23} and \eqref{eq:32} can be rewritten as, respectively,
   \begin{equation}
   \label{eq:23 r=4}
   a(p)_{23} ~=~ \frac{p^3(p^4 - 1)}{120} a(2)_{23}
   \end{equation}
and
   \begin{equation}
   \label{eq:32 r=4}
   \begin{split}
   \frac{h}{D}\bigl(a(p)_{22} & ~-~ a(p)_{33}\bigr)  ~+~ \frac{240}{D}a(p)_{32} \\
   &~=~ \frac{h}{D}\cdot \frac{p^3(p^4 - 1)}{120}\bigl(a(2)_{22} ~-~ a(2)_{33}\bigr) ~+~ 2\cdot \frac{p^3(p^4 - 1)}{D}a(2)_{32},
   \end{split}
   \end{equation}
where $D = D(4, h) = \gcd(h, 240)$.  Since $h/D$ and $240/D$ are relatively prime, we can make the same argument as in the case $r = 2$, $h = 5$.  More precisely, the two conditions, \eqref{eq:23 r=4} and \eqref{eq:32 r=4}, above tell us that to choose an element of $\Derbar(R)$, one first chooses five arbitrary elements in $2\bZ$ for the $(1, 1)$, $(2, 2)$, $(2, 3)$, $(3, 2)$, and $(3, 3)$ entries in $f(2)$.  Then for each odd prime $p$, one chooses an element in $p \bZ$ for the $(1, 1)$ entry in $f(p)$ and two more arbitrary integers to form $f(p)$.  These last two degrees of freedom come from \eqref{eq:32 r=4} and are completely analogous to \eqref{eq:freedom}.  This allows us to make the identification \eqref{eq:der r=2 h=5} for $\Derbar(R) = \Derbar(S((p^4), h)$.

To complete the proof, observe that the $p = 2$ component of a typical element  in $\Innbar(R)$ (see \eqref{eq:inn}) has the form
   \[
   \begin{bmatrix}
   0 & 0 & 0 \\
   0 & -sh & -240s \\
   0 & h(j - t) + 240k & sh
   \end{bmatrix}.
   \]
Moreover, we have that
   \[
   \lbrace sh \colon s \in 2\bZ \rbrace ~=~ 2h \bZ
   \]
and 
   \[
   \lbrace h(j - t) + 240k \colon j \equiv t \pmod{2} \rbrace ~=~ 2D\bZ.
   \]
Therefore, it follows as in the case when $r = 2$ and $h = 5$ that there is a group isomorphism
   \[
   \Hl^1(S((p^4),\, h)) ~\cong~ \frac{2\bZ}{2h\bZ} \times \frac{2\bZ}{2D\bZ} \times \bZ^{\bN},    
   \]
as desired.

This proves Theorem \ref{thm:H1 n=3} when $r = 4$, except for the commutativity assertion, to which we now turn.


\subsection{Proof of the commutativity assertion in Theorem \ref{thm:H1 n=3}}
\label{subsec:commutativity}

We can now show that the graded algebra $\Hl^{\leq 1}(R) = \Hl^{\leq 1}(S((p^r),\, h))$ is commutative if and only if $D = 1$, i.e.\ $h$ and $2^r(2^r - 1)$ are relatively prime.  Since Theorem \ref{thm:H0 n=3} already established the commutativity of $\Hl^0(R)$, commutativity of $\Hl^{\leq 1}(R)$ here means that elements in $\Hl^0(R)$ commute with those in $\Hl^1(R)$.

To do this, let $g \in \Hl^0(R) \subseteq \Endbar(R) = \cF^0(R)$ be a $0$-cocycle, as in Theorem \ref{thm:H0 n=3} (2), and let $f$ $=$ $\lbrace f(p) = (f(p)_{ij}) \rbrace_{p \in \cP}$ be a $1$-cocycle, i.e.\ an element of $\Derbar(R)$.  Then the map 
   \[
   (g \circ f ~-~ f \circ g)(p) ~=~ gf(p) - f(p)g
   \]
is represented by the matrix
   \begin{equation}
   \label{eq:diff}
   \begin{bmatrix}
   0 & 0 & 0 \\
   0 & -kf(p)_{23} & (j - t)f(p)_{23} \\
   0 & kf(p)_{22} - kf(p)_{33} + (t - j)f(p)_{32}  & kf(p)_{23}
   \end{bmatrix}
   \end{equation}
Notice that in all $64$ cases under consideration, we have
   \begin{equation}
   \label{eq:G}
   \begin{split}
   G &~=~ \gcd(b_p(b_p - 1))_{p \in \cP} ~=~ \gcd(p^r(p^r - 1))_{p \in \cP} \\
      &~=~ 2^r(2^r - 1) \\
      &~=~ \begin{cases} 2 & \text{if } r = 1\\
                                       12 & \text{if } r = 2 \\
                                       240 & \text{if } r = 4. \end{cases}
   \end{split}
   \end{equation}
In particular, $D = D(r, h) = 1$ if and only if $h$ and $G$ are relatively prime.

Suppose that $D = 1$.  We need to show that the cohomology classes represented by the $1$-cocycles $g \circ f$ and $f \circ g$ are the same, i.e.\ the $1$-cocycle $(g\circ f - f\circ g)$ is actually a $1$-coboundary.  Since the entries in $g$ satisfies
   \begin{equation}
   \label{eq:hG}
   h(t - j) ~=~ kG,
   \end{equation}
it follows that $k \equiv 0 \pmod{h}$.  Using Lemma \ref{lem1:H0 n=3}, one observes that the matrix 
   \[
   \beta ~=~
    \frac{k}{h}
   \begin{bmatrix}
   0 & 0 & 0 \\
   0 & f(2)_{22} & f(2)_{23} \\
   0 & f(2)_{32} & f(2)_{33}
   \end{bmatrix}
   \]
represents a $0$-cochain in $\cF^*(R)$.  We claim that
   \begin{equation}
   \label{eq:beta}
   d^0 \beta ~=~ g\circ f ~-~ f\circ g.
   \end{equation}
To see this, first note that the matrix in \eqref{eq:inn} is the $p$th component of a $1$-coboundary in $\cF^*(R)$.   Applying this to the $0$-cochain $\beta$,  we see that the first column and first row of $(d^0\beta)(p)$ are both $0$ and that its lower-right $2$-by-$2$ submatrix is
   \begin{equation}
   \label{eq:dbeta}
   \frac{k}{h}
   \begin{bmatrix}
   -f(2)_{23}c_p & -f(2)_{23}p^r(p^r - 1) \\
   f(2)_{32}p^r(p^r - 1) + (f(2)_{22} - f(2)_{33})c_p & f(2)_{23}c_p
   \end{bmatrix}.
   \end{equation}
To prove \eqref{eq:beta}, we only need to show that the matrix in \eqref{eq:dbeta} coincides with the lower-right $2$-by-$2$ block of the matrix \eqref{eq:diff}.  Now \eqref{eq:23} and \eqref{eq:G} imply that
   \begin{equation}
   \label{eq:fp23}
   f(p)_{23} ~=~ \frac{p^r(p^r - 1)}{G}f(2)_{23}.
   \end{equation}
One infers from \eqref{eq:fp23} and \eqref{eq:hG} that
   \[
   \begin{split}
   -\frac{k}{h}f(2)_{23}p^r(p^r - 1) 
   &~=~ -\frac{k}{h}Gf(p)_{23} \\
   &~=~ (j - t)f(p)_{23}.
   \end{split}
   \]
This shows that the $(2, 3)$ entries in $d^0 \beta$ and $(g\circ f - f\circ g)$ coincide.

Similarly, we have
   \[
   \begin{split}
   \frac{k}{h}f(2)_{23}c_p 
   &~=~ \frac{k}{h}f(2)_{23} \frac{hp^r(p^r - 1)}{G} \\
   &~=~ kf(p)_{23}.
   \end{split}
   \]
This shows that the $(3, 3)$ (and hence also the $(2, 2)$) entries in $d^0 \beta$ and $(g\circ f - f\circ g)$ coincide.

Finally, we have
   \[
   \begin{split}
   \frac{k}{h} &\Bigl(c_p(f(2)_{22} - f(2)_{33})
   ~+~ p^r(p^r - 1)f(2)_{32} \Bigr) \\
   &~=~ \frac{k}{h} \Bigl(c_2(f(p)_{22} - f(p)_{33}) ~+~ 2^r(2^r - 1)f(p)_{32} \Bigr) \quad \text{by } \eqref{eq:32} \\
   &~=~ \frac{k}{h} \Bigl(h(f(p)_{22} - f(p)_{33}) ~+~ Gf(p)_{32}\Bigr) \\
   &~=~ k(f(p)_{22} - f(p)_{33}) ~+~ (t - j)f(p)_{32}.
   \end{split}
   \]
The second equality follows from the fact that $c_2 = h$ and $G = 2^r(2^r - 1)$.  This shows that the $(3, 2)$ entries in $d^0 \beta$ and $(g\circ f - f\circ g)$ coincide, and thus the claim \eqref{eq:beta} is proved.  Of course, this implies that the cohomology classes represented by the $1$-cocycles $g\circ f$ and $f\circ g$ are equal.

Now suppose that $D = \gcd(h,\, 2^r(2^r - 1)) > 1$.  (This can only happen when $r \not= 1$.)  We must show that the graded algebra $\Hl^{\leq 1}(R)$ is not commutative.  Consider the $0$-cocycle
   \[
   g ~=~ \begin{bmatrix} 0 & 0 & 0 \\
                         0 & 0 & 0 \\
                         0 & \frac{h}{D} & \frac{G}{D} \end{bmatrix}.
   \]
That this is indeed a $0$-cocycle follows from Theorem \ref{thm:H0 n=3} (2).  Also, by Lemma \ref{lem1:H1 n=3}, the sequence (indexed by the primes) of linear self maps on $R$ represented by the matrices 
   \[
   f ~=~ \left\lbrace \begin{bmatrix}
         0 & 0 & 0 \\
         0 & 0 & \frac{2p^r(p^r - 1)}{G} \\
         0 & 0 & 0 \end{bmatrix}
         \right\rbrace_{p \in \cP}
   \]
is an element of $\Derbar(R)$, i.e.\ a $1$-cocycle in $\cF^1(R)$.  We claim that the $1$-cocycle $(g \circ f - f \circ g)$ is not a $1$-coboundary, which would imply that the graded algebra $\Hl^{\leq 1}(R)$ is not commutative.  To see this, observe from \eqref{eq:diff} that the $(3,\, 3)$ entry in the matrix of $(gf(2) - f(2)g)$ is
   \[
   \frac{h}{D} \cdot \frac{2 \cdot 2^r(2^r - 1)}{G} ~=~ \frac{2h}{D}.
   \]
 From Lemma \ref{lem2:H1 n=3}, the $(3, 3)$ entry of the component for $p = 2$ of an element in $\Innbar(R)$ is of the form  $sc_2 = sh$ for some even integer $s$.  But then the equality
   \[
   sh ~=~ \frac{2h}{D}
   \]
would imply 
   \[
   Ds ~=~ 2,
   \]
which is absurd, since $D > 1$.  In other words, we have shown that the cohomology classes $g \in \Hl^0(R)$ and $\lbrack f \rbrack \in \Hl^1(R)$ do not commute with each other, as desired.

This finishes the proof of the commutativity assertion in Theorem \ref{thm:H1 n=3}.  The proof of Theorem \ref{thm:H1 n=3} is now complete.


\section{$\Hl^0$ of certain filtered $\lambda$-ring structures on $\bZ \lbrack x \rbrack/(x^4)$}
\label{sec:n=4}

The purpose of this section is to compute the algebras $\Hl^0(R)$ for the following $61$ filtered $\lambda$-ring structures $R$ on the truncated polynomial algebra $\bZ \lbrack x \rbrack/(x^4)$, with $x$ in a fixed positive filtration:

\begin{enumerate}
\item $R = K(\mathbf{CP}^3)$ with Adams operations
       \[ 
       \psi^p(x) ~=~ (1 + x)^p ~-~ 1
       \]
for $p$ primes.

\item $R = S(h, d_2)$ with
      \[
      S(h,\, d_2) ~=~ \left\lbrace \psi^p(x) ~=~ p^2x ~+~ h \frac{p^2(p^2 - 1)}{12}x^2 ~+~ d_px^3 \right\rbrace_{p \in \cP},
      \]
where $h \in \lbrace 1, 5 \rbrace$, $d_2 \in \lbrace 0, 2, 4, \ldots , 58\rbrace$, and
      \begin{equation}
      \label{eq:d}
      d_p ~=~ \frac{p^2(p^4 - 1)}{60}d_2 ~+~ \frac{p^2(p^2 - 1)(p^2 - 4)}{360}h^2
      \end{equation}
for odd primes $p$.
\end{enumerate}
In this notation, the $K$-theory filtered $\lambda$-ring of $\mathbf{HP}^3$, the quaternionic projective $3$-space, is $S(1, 0)$.

It is shown in \cite{yau4} that the truncated polynomial algebra $\bZ \lbrack x \rbrack/(x^4)$ admits uncountably many isomorphism classes of filtered $\lambda$-ring structures.   Among those classes are the $61$ filtered $\lambda$-rings above.  Moreover, if $X$ is a torsionfree space (i.e.\ its integral cohomology is $\bZ$-torsionfree) whose $K$-theory filtered ring is the truncated polynomial ring $\bZ \lbrack x \rbrack/(x^4)$, then $K(X)$ is isomorphic as a filtered $\lambda$-ring to one of the above $61$ filtered $\lambda$-rings.  In other words, among the uncountably many isomorphism classes of filtered $\lambda$-ring structures on $\bZ \lbrack x \rbrack/(x^4)$, at most $61$ of them, those listed above, can possibly be topologically realized by torsionfree spaces.

As in previous sections, we will use the standard $\bZ$-basis $\lbrace 1, x, x^2, x^3 \rbrace$ for the ring $\bZ \lbrack x \rbrack/(x^4)$, and, using this basis, each $\bZ$-linear self map is represented by an element in $M(4,\bZ)$, the algebra of $4$-by-$4$ matrices with integer entries.  The algebra $\Hl^0(R)$ then becomes a subalgebra of the matrix algebra $M(4,\bZ)$.

Here are the main results of this section.

\begin{thm}
\label{thm1:H0 n=4}
With the standard basis $\lbrace 1, x, x^2, x^3 \rbrace$ of $\bZ \lbrack x \rbrack/(x^4)$, $\Hl^0(K(\mathbf{CP}^3))$ is the subalgebra of $M(4, \bZ)$ consisting of matrices of the form
   \begin{equation}
   \label{eq1:H0 n=4}
   \begin{bmatrix}
   a & & & \\
     & j & & \\
     & k & j + 2k & \\
     & l & 4k + 6l & j + 6k + 6l
   \end{bmatrix}.
   \end{equation}
This algebra is commutative and, as an additive group, is free of rank $4$ over $\bZ$.
\end{thm}

\begin{thm}
\label{thm2:H0 n=4}
Consider any one of the $60$ filtered $\lambda$-rings $S(h, d_2)$.  Then $\Hl^0(S(h, d_2))$ is the subalgebra of $M(4, \bZ)$ consisting of matrices of the form
   \begin{equation}
   \label{eq2:H0}
   \begin{bmatrix}
   a & & & \\
     & j & & \\
     & k & r & \\ 
     & l & s & w 
   \end{bmatrix}
   \end{equation}
such that the following is true:
   \begin{enumerate}
   \item If $h = 1$, then
   \begin{equation}
   \label{eq:rsw}
   \begin{split}
   r &~=~ j + 12k, \\
   (6d_2 + 1)s &~=~ (8 - 12d_2)k ~+~ 60l,\, \text{and} \\
   w &~=~ 6s ~+~ j ~+~ 12k.
   \end{split}
   \end{equation}
   \item If $h = 5$, then 
   \begin{equation}
   \label{eq:krsw}
   \begin{split}
   k &~\equiv~ 0 \pmod{5}, \\
   r &~=~ j + \frac{12k}{5}, \\
   (6d_2 + 25)s &~=~ (200 - 12d_2)k ~+~ 300l,\, \text{and} \\
   (6d_2 + 25)w &~=~ (6d_2 + 25)j ~+~ 300k ~+~ 360l.
   \end{split}
   \end{equation}
   \end{enumerate}
In each of the $60$ cases, the algebra $\Hl^0(S(h, d_2))$ is commutative and, as an additive group, is free of rank $4$ over $\bZ$.
\end{thm}

To prove these results, we begin by computing the group $\cF^0(R)$ of $0$-cochains.

\begin{lemma}
\label{lem:H0 n=4}
Let $R$ be any filtered $\lambda$-ring structure on $\bZ \lbrack x \rbrack/(x^4)$.  Then the group $\Endbar(R) = \cF^0(R)$ consists of the matrices
   \begin{equation}
   \label{eq:endbar n=4}
   g ~=~ 
   \begin{bmatrix}
   a & & & \\
     & j & n & u \\
     & k & r & v \\
     & l & s & w 
   \end{bmatrix}
   \end{equation}
such that
   \begin{enumerate}
   \item $n \equiv u \equiv 0 \pmod{6}$,
   \item $r - j \equiv s \equiv 0 \pmod{2}$,
   \item $w - j \equiv v \equiv 0 \pmod{3}$.
   \end{enumerate} 
\end{lemma}

\begin{proof}
Let $g$ be a $\bZ$-linear self map on $R$ represented by the matrix
   \[
   g ~=~ 
   \begin{bmatrix}
   a & i & m & t \\
   b & j & n & u \\
   c & k & r & v \\
   d & l & s & w
   \end{bmatrix}.
   \]
By definition, $g \in \Endbar(R)$ if and only if it satisfies \eqref{eq:endbar}, which is equivalent to \eqref{eq:endbar p} for $0 \leq i \leq 3$.  For $i = 0$, \eqref{eq:endbar p} says 
   \[
   \begin{split}
   g(1^p) 
   &~=~ g(1) ~=~ a + bx + cx^2 + dx^3 \\
   &~\equiv~ (a + bx + cx^2 + dx^3)^p \pmod{p} \\
   &~\equiv~ a + bx^p \pmod{p} \\
   &~=~ a \pmod{p} \quad \text{ if } p \geq 5.
   \end{split}
   \]
This implies that 
   \[
   b ~=~ c ~=~ d ~=~ 0;
   \]
that is, $g(1) = a \in \bZ$.  The cases $p = 2$ and $3$ give no additional information.

When $i = 1$, the condition \eqref{eq:endbar p} says
   \[
   \begin{split}
   g(x)^p 
   &~=~ (i + jx + kx^2 + lx^3)^p \\
   &~\equiv~ i \pmod{p}  \quad \text{ if } p \geq 5\\
   &~\equiv~ g(x^p) \pmod{p} \\ 
   &~=~ 0 \pmod{p} \quad \text{ if } p \geq 5.
   \end{split}
   \]
This implies that $i = 0$.  We have not used the cases $p = 2$ and $3$ yet.  If $p = 2$, then 
   \[
   \begin{split}
   g(x)^2 
   &~\equiv~ jx^2 \pmod{2} \\
   &~\equiv~ g(x^2) \pmod{2} \\
   &~=~ m + nx + rx^2 + sx^3.
   \end{split}
   \]
This implies that
   \[
   m ~\equiv~ n ~\equiv~ s ~\equiv~ j - r ~\equiv~ 0 \pmod{2}.
   \]
If $p = 3$, then
   \[
   \begin{split}
   g(x)^3
   &~\equiv~ jx^3 \pmod{3} \\
   &~\equiv~ g(x^3) \pmod{3} \\
   &~=~ t + ux + vx^2 + wx^3.
   \end{split}
   \]
This implies that
   \[
   t ~\equiv~ u ~\equiv~ v ~\equiv~ j - w ~\equiv~ 0 \pmod{3}.
   \]

For $i = 2$, the condition \eqref{eq:endbar p} says
   \[
   \begin{split}
   g(x^2)^p 
   &~=~ (m + nx + rx^2 + sx^3)^p \\
   &~\equiv~ m + nx^p \pmod{p} \\
   &~\equiv~ g(x^{2p}) \pmod{p} \\
   &~=~ 0.
   \end{split}
   \]
This implies that $m = 0$ and $n \equiv 0 \pmod{6}$.

Finally, for $i = 3$, \eqref{eq:endbar p} says
   \[
   \begin{split}
   g(x^3)^p
   &~=~ (t + ux + vx^2 + wx^3)^p \\
   &~\equiv~ t + ux^p \pmod{p} \\
   &~\equiv~ g(x^{3p}) \pmod{p} \\
   &~=~ 0.
   \end{split}
   \]
So $t = 0$ and $u \equiv 0 \pmod{6}$.

This proves the Lemma.
\end{proof}

\begin{proof}[Proof of Theorem \ref{thm1:H0 n=4}]
Let $R$ denote $K(\mathbf{CP}^3)$ and let $g \in \Endbar(R)$ be as in \eqref{eq:endbar n=4}.  Then $g \in \Hl^0(R)$ if and only if $g\psi^p = \psi^p g$ for all primes $p$.  Since the matrices for both $g$ and $\psi^p$ have $0$'s in their first rows and first columns, except in the $(1, 1)$ entries, $g$ commutes with $\psi^p$ if and only if their lower-right $3$-by-$3$ submatrices commute.  Note that the matrix for $\psi^p$ is 
   \[
   \psi^p ~=~
   \begin{bmatrix}
   1 & & & \\
     & p & & \\
     & \binom{p}{2} & p^2 & \\
     & \binom{p}{3} & p^2(p - 1) & p^3
   \end{bmatrix}.
   \]
Denote by $A$ (respectively $B_p$) the lower-right $3$-by-$3$ submatrix of $g$ (respectively $\psi^p$).  Then $A$ commutes with $B_p$ means that
\begin{small}
   \begin{equation}
   \label{eq:abp}
   \begin{split}
   & AB_p ~=~
   \begin{bmatrix}
   pj + \binom{p}{2}n + \binom{p}{3}u & p^2n + p^2(p - 1)u & p^3u \\
   pk + \binom{p}{2}r + \binom{p}{3}v & p^2r + p^2(p - 1)v & p^3v \\
   pl + \binom{p}{2}s + \binom{p}{3}w & p^2s + p^2(p - 1)w & p^3w  
   \end{bmatrix} \\
   &~=~ 
   \begin{bmatrix}
   pj & pn & pu \\
   \binom{p}{2}j + p^2k & \binom{p}{2}n + p^2r & \binom{p}{2}u + p^2v \\
   \binom{p}{3}j + p^2(p - 1)k + p^3l & \binom{p}{3}n + p^2(p - 1)r + p^3s & \binom{p}{3}u + p^2(p - 1)v + p^3w
   \end{bmatrix} \\
   &~=~ B_pA.
   \end{split}
   \end{equation} 
\end{small}

For any prime $p$, comparing the $(1, 3)$ entries in \eqref{eq:abp} yields $u = 0$.  Applying this to the $(1, 2)$ and the $(2, 3)$ entries gives $n = v = 0$.  In particular, if $g$ commutes with $\psi^p$ for any $p$, then $g$ must be lower-triangular, i.e.\ it respects the filtration of $R$.  The diagonal entries in \eqref{eq:abp} provide no new information about the entries of $g$.  Since $v = 0$, the $(2, 1)$ entries in \eqref{eq:abp} yields 
   \begin{equation}
   \label{eq:r}
   r ~=~ j + 2k,
   \end{equation}
as stated in the Theorem.

Now consider the $(3, 1)$ entries in \eqref{eq:abp}.  If $p = 2$, then, since $\binom{2}{3} = 0$, one obtains 
   \begin{equation}
   \label{eq:s}
   s ~=~ 4k + 6l.
   \end{equation}
If $p = 3$, then one obtains the equation
   \[
   3l + 3s + w ~=~ j + 18k + 27l,
   \]
which is equivalent to
   \begin{equation}
   \label{eq:w}
   w ~=~ j + 6k + 6l.
   \end{equation}
The $(3, 2)$ entries in \eqref{eq:abp} give no new information, since the equation obtained from them is
   \[
   w ~=~ r + s
   \]
for each prime $p$.  It is clear from \eqref{eq:r}, \eqref{eq:s}, and \eqref{eq:w} that this is true.  Conversely, it is straightforward to check that \eqref{eq:r}, \eqref{eq:s}, and \eqref{eq:w} imply that $AB_p$ and $B_pA$ have the same $(3, 1)$ and $(3, 2)$ entries for all primes $p$.

This proves that $\Hl^0(R)$ is the subalgebra of $M(4, \bZ)$ consisting of matrices of the form \eqref{eq1:H0 n=4}.  It is now also easy to see that, as an additive group, $\Hl^0(R)$ is free of rank $4$ over $\bZ$, since $a$, $j$, $k$, and $l$ are arbitrary and the other three entries are linear combinations in $j$, $k$, and $l$.

To show that $\Hl^0(R)$ is a commutative algebra, let $g^\prime$ be another element in it of the form \eqref{eq1:H0 n=4} with $a^\prime$, $j^\prime$, $k^\prime$, and $l^\prime$ in place of $a$, $j$, $k$, and $l$, respectively.  It is clearly enough to show that the $(i, j)$ entries in $gg^\prime$ and $g^\prime g$ coincide for $(i, j) = (3, 2)$, $(4, 2)$, and $(4, 3)$.  Denote the $(i, j)$ entry in a matrix $A$ by $A_{ij}$.  A little bit of matrix computation then shows that $(gg^\prime)_{ij}$ is
   \begin{itemize}
   \item $kj^\prime + k^\prime j + 2kk^\prime$ if $(i, j) = (3, 2)$, 
   \item $lj^\prime + l^\prime j + 4kk^\prime + 6(k^\prime l + l^\prime k + ll^\prime)$ if $(i, j) = (4, 2)$, and
   \item $(4k + 6l)(j^\prime + 2k^\prime) + (4k^\prime + 6l^\prime)(j + 2k) + (4k + 6l)(4k^\prime + 6l^\prime)$ if $(i, j ) = (4, 3)$.
   \end{itemize}
Since each of these entries remains the same if $j$ (respectively $k$, $l$) and $j^\prime$ (respectively, $k^\prime$ and $l^\prime$) are interchanged, it follows that $gg^\prime$ is equal to $g^\prime g$, as desired.

This finishes the proof of Theorem \ref{thm1:H0 n=4}.
\end{proof}

\begin{proof}[Proof of Theorem \ref{thm2:H0 n=4}]
Let $R$ be $S(h, d_2)$.  This proof is quite similar to the proof of Theorem \ref{thm1:H0 n=4}.

The Adams operation $\psi^p$ in $R$ is given by
   \[
   \psi^p(x) ~=~ p^2x + h\frac{p^2(p^2 - 1)}{12}x^2 + d_px^3
   \]
and $\psi^p$ is a ring map on $R$.  Therefore, the matrix of $\psi^p$ has $0$'s in its first row and first column, except for the entry $1$ in the $(1,1)$ spot, and its lower-right $3$-by-$3$ submatrix is
   \[
   B_p ~=~ \begin{bmatrix}
   p^2 & 0 & 0 \\
   hp^2(p^2 - 1)/12 & p^4 & 0 \\
   d_p & hp^4(p^2 - 1)/6 & p^6
   \end{bmatrix}.
   \]
Let $g \in \Endbar(R)$ be a $0$-cochain whose matrix is as in \eqref{eq:endbar n=4}, and let $A$ denote its lower-right $3$-by-$3$ submatrix.  Then $g$ commutes with $\psi^p$ if and only if $A$ commutes with $B_p$.  As in \eqref{eq:abp}, this means that
   \begin{equation}
   \label{eq2:abp}
   \begin{split}
   &
   \begin{bmatrix}
   p^2j + hp^2(p^2 - 1)n/12 + d_pu & p^4n + hp^4(p^2 - 1)u/6 & p^6u \\
   p^2k + hp^2(p^2 - 1)r/12 + d_pv & p^4r + hp^4(p^2 - 1)v/6 & p^6v \\
   p^2l + hp^2(p^2 - 1)s/12 + d_pw & p^4s + hp^4(p^2 - 1)w/6 & p^6w \\
   \end{bmatrix} ~=~ \\
   &
   \begin{bmatrix}
   p^2j & p^2n & p^2u \\
   hp^2(p^2 - 1)j/12 + p^4k & hp^2(p^2 - 1)n/12 + p^4r & hp^2(p^2 - 1)u/12 + p^4v \\
   d_pj + \frac{hp^4(p^2 - 1)k}{6} + p^6l & d_pn + \frac{hp^4(p^2 - 1)r}{6} + p^6s & d_pu + \frac{hp^4(p^2 - 1)v}{6} + p^6w
   \end{bmatrix}
   \end{split}
   \end{equation}

Again, for any prime $p$, the $(1, 3)$ entries tell us that $u = 0$.  Applying this to the $(1, 2)$ and the $(2, 3)$ entries, one obtains $n = v = 0$, i.e.\ $g$ is lower-triangular.  The diagonal entries in \eqref{eq2:abp} give no new information about the entries in $g$.

The $(2, 1)$ entries in \eqref{eq2:abp} yield, for any prime $p$, the equation
   \[
   h(r - j) ~=~ 12k.
   \]
If $h = 1$, this means that
   \[
   r ~=~ j + 12k,
   \]
as stated in the Theorem.  If $h = 5$, then $5$ must divide $k$ and
   \[
   r ~=~ j + \frac{12k}{5},
   \]
as desired.

Setting $p = 2$, the $(3, 1)$ and $(3, 2)$ entries in \eqref{eq2:abp} yield the simultaneous equations
   \begin{equation}
   \label{eq:sh}
   \begin{split}
   hs + d_2 w &~=~ d_2j + 8hk + 60l \\
   6s - hw &~=~ -hj - 12k
   \end{split}
   \end{equation}
in $s$ and $w$.  It is an elementary exercise to see that this system of linear equations gives rise to the required conditions, \eqref{eq:rsw} and \eqref{eq:krsw}, for $s$ and $w$ in the statement of Theorem \ref{thm2:H0 n=4}.  Furthermore, using \eqref{eq:d}, which expresses each $d_p$ in terms of $d_2$ and $h$, a slightly tedious but easy calculation shows that \eqref{eq:rsw} and \eqref{eq:krsw} make the $(3, 1)$ (respectively $(3, 2)$) entries in $AB_p$ and $B_pA$ coincide for any prime $p$.  This proves that $\Hl^0(R)$ consists of exactly the matrices \eqref{eq2:H0} satisfying \eqref{eq:rsw} or \eqref{eq:krsw}, depending on whether $h$ is $1$ or $5$.

To see that $\Hl^0(R)$ is a commutative algebra, note that
   \[
   \begin{split}
   r &~=~ j + a_1 k \\
   s &~=~ a_2k + a_3l \\
   w &~=~ j + a_4k + a_5l
   \end{split}
   \]
for some rational numbers $a_1, \ldots , a_5$.  Therefore, the proof for the commutativity of $\Hl^0(K(\mathbf{CP}^3))$ can be used virtually verbatim here as well.

To see that $\Hl^0(R)$ is free of rank $4$, first consider the case $h = 1$.  We will use the notations in \eqref{eq:rsw}.  Since $r$ (respective $w$) are linear combinations in $j$ and $k$ (respectively $s$, $j$, and $k$), it suffices to show that 
   \begin{equation}
   \label{eq:skl}
   \lbrace (k, l, s) \in \bZ^{\times 3} \colon  (6d_2 + 1)s ~=~ (8 - 12d_2)k + 60l \rbrace
   \end{equation}
is free of rank $2$.  Define $\alpha$ by
   \[
   \alpha ~=~ \gcd(8 - 12d_2,\, 60,\, 6d_2 + 1).
   \]
Then there is a surjective $\bZ$-linear map
   \[
   \varphi \colon \bZ^{\times 3} ~\to~ \alpha \bZ ~\cong~ \bZ,
   \]
where
   \[
   \varphi((k, l, s)) ~=~ (8 - 12d_2)k + 60l - (6d_2 + 1)s.
   \]
The kernel of $\varphi$ is exactly the group in \eqref{eq:skl}, which shows that it has rank $2$, as desired.

Now consider the case $h = 5$.  Let $H$ be the subgroup of 
   \[
   K ~=~ \bZ \times (5\bZ) \times \bZ^{\times 3} ~\cong~ \bZ^{\times 5} 
   \]
consisting of elements $(j, k, l, s, w) \in K$ for which the last two equations in \eqref{eq:krsw} hold.  As in the case $h = 1$, it suffices to show that $H$ is free of rank $3$.  Define $\beta$ and $\gamma$ by
   \[
   \begin{split}
   \beta &~=~ \gcd(5(200 - 12d_2),\, 300,\, 6d_2 + 25) \\
   \gamma &~=~ \gcd(5(300),\, 360,\, 6d_2 + 25).
   \end{split}
   \]
Then there is a surjective $\bZ$-linear map
   \[
   \phi \colon K ~\to~ \beta \bZ \times \gamma \bZ ~\cong~ \bZ^{\times 2},
   \]
where
   \[
   \begin{split}
   \phi((j, k, l, s, w)) ~=~ & \bigl( (200 - 12d_2)k + 300l - (6d_2 + 25)s, \\
   & (6d_2 + 25)(j - w) + 300k + 360l \bigr).
   \end{split}
   \]
The kernel of $\phi$ is exactly $H$, which shows that $H$ is free of rank $3$, as desired.

This finishes the proof of Theorem \ref{thm2:H0 n=4}.
\end{proof}




\begin{thebibliography}{99}

\bibitem{at}M. F. Atiyah and D. O. Tall, Group representations, $\lambda$-rings and the $J$-homomorphism, Topology 8 (1969), 253-297.


\bibitem{ger}M.\ Gerstenhaber, On the deformation of rings and algebras, Ann.\ Math.\ (2) 79 (1964), 59-103.


\bibitem{gro}A.\ Grothendieck, Classes de faisceaux et th\'eor\`eme de Riemann-Roch, in: SGA 6, 22-70, Lecture Notes in Math., vol. 225, Springer, Berlin, 1971.

\bibitem{hazewinkel}M. Hazewinkel, Formal groups and applications,  Pure and Applied Mathematics 78, Academic Press, New York-London, 1978. 

\bibitem{knutson}D. Knutson, $\lambda $-rings and the representation theory of the symmetric group, Lecture Notes in Math. 308, Springer-Verlag, Berlin-New York, 1973.

\bibitem{patras}F. Patras, Lambda-rings, in: Handbook of algebra, vol. 3, 961-986, 
North-Holland, Amsterdam, 2003. 


\bibitem{wil}C. Wilkerson, Lambda-rings, binomial domains, and vector bundles over $CP(\infty)$, Comm. Algebra 10 (1982), 311-328.

\bibitem{yau1}D. Yau, Moduli space of filtered $\lambda$-ring structures over a filtered ring, Int. J. Math. Math. Sci. (2004), no. 39, 2065-2084.


\bibitem{yau3}D. Yau, Cohomology of $\lambda$-rings, J.\ Algebra, accepted for publication.

\bibitem{yau4}D. Yau, On $\lambda$-rings and topological realization, preprint, available at \texttt{http://www.math.uiuc.edu/$\sim$dyau/papers.html}.

\end{thebibliography}
\end{document}